\documentclass[12pt]{compositio}

\usepackage{amssymb}
 \usepackage{amsmath}
\usepackage[all]{xy}
\usepackage{eucal}
\usepackage{xfrac}
\usepackage{relsize}
\usepackage{stmaryrd}

\usepackage{tikz}
\usepackage{tikz-cd}
\usetikzlibrary{calc,matrix}
\tikzset{w/.style={circle, draw,inner sep=1pt},b/.style={circle,draw,fill,inner sep=2pt}, s/.style={rectangle, draw,inner sep=3pt}}

\usepackage[colorlinks=true]{hyperref}
\usepackage{cleveref}

\newcounter{conj}
\theoremstyle{plain}
\newtheorem{thm}{Theorem}[section]
\crefname{thm}{Theorem}{Theorem}
\newtheorem{prop}[thm]{Proposition}
\crefname{prop}{Proposition}{Proposition}
\newtheorem{lm}[thm]{Lemma}
\crefname{lm}{Lemma}{Lemma}
\newtheorem{cor}[thm]{Corollary}
\crefname{cor}{Corollary}{Corollary}
\newtheorem*{theorem}{Theorem}
\newtheorem*{corollary}{Corollary}

\theoremstyle{definition}
\newtheorem{conjecture}[conj]{Conjecture}
\crefname{conjecture}{Conjecture}{Conjecture}
\newtheorem{defn}[thm]{Definition}
\crefname{defn}{Definition}{Definition}

\theoremstyle{remark}
\newtheorem{remark}[thm]{Remark}
\crefname{remark}{Remark}{Remark}
\newtheorem{warn}[thm]{Warning}
\crefname{warn}{Warning}{Warning}
\newtheorem{example}[thm]{Example}
\crefname{example}{Example}{Example}

\newcommand{\bP}{\mathbb{P}}
\newcommand{\bG}{\mathbb{G}}
\newcommand{\bA}{\mathbb{A}}

\newcommand{\Z}{\mathbf{Z}}
\newcommand{\C}{\mathbf{C}}

\newcommand{\rZ}{\mathrm{Z}}
\newcommand{\Sp}{\mathcal{S}}

\renewcommand{\d}{\mathrm{d}}
\newcommand{\dR}{\mathrm{dR}}
\newcommand{\DR}{\mathrm{DR}}
\newcommand{\CE}{\mathrm{CE}}
\newcommand{\cCE}{\widehat{\mathrm{CE}}}
\newcommand{\aug}{\mathrm{aug}}
\newcommand{\adic}{\mathrm{adic}}
\newcommand{\dPSt}{\mathrm{dPSt}}

\newcommand{\Thick}{\mathrm{Thick}}

\newcommand{\fib}{\mathrm{fib}}
\newcommand{\cofib}{\mathrm{cofib}}
\newcommand{\Fun}{\mathrm{Fun}}
\newcommand{\id}{\mathrm{id}}
\newcommand{\ev}{\mathrm{ev}}
\newcommand{\LieAlgd}{\mathrm{LieAlgd}}
\newcommand{\QCoh}{\mathrm{QCoh}}
\newcommand{\Coh}{\mathrm{Coh}}
\newcommand{\Tot}{\mathrm{Tot}}
\newcommand{\Perf}{\mathrm{Perf}}
\newcommand{\APerf}{\mathrm{APerf}}
\newcommand{\Symp}{\mathrm{Symp}}
\newcommand{\Lag}{\mathrm{Lag}}
\newcommand{\Pois}{\mathrm{Pois}}
\newcommand{\Cois}{\mathrm{Cois}}
\newcommand{\TCois}{\mathrm{TCois}}
\newcommand{\Comp}{\mathrm{Comp}}
\newcommand{\Obs}{\mathrm{Obs}}
\newcommand{\Isot}{\mathrm{Isot}}

\newcommand{\Sym}{\mathrm{Sym}}
\newcommand{\Pol}{\mathrm{Pol}}
\newcommand{\wPol}{\mathbf{Pol}}
\newcommand{\LagThick}{\mathrm{LagThick}}
\newcommand{\Hom}{\mathrm{Hom}}
\newcommand{\Com}{\mathrm{Com}}
\newcommand{\Lie}{\mathrm{Lie}}
\newcommand{\Fol}{\mathrm{Fol}}
\newcommand{\Op}{\mathrm{Op}}
\newcommand{\preLie}{\mathrm{preLie}}
\newcommand{\Conv}{\mathrm{Conv}}
\newcommand{\End}{\mathrm{End}}
\newcommand{\co}{\mathrm{co}}
\newcommand{\cu}{\mathrm{cu}}
\renewcommand{\nu}{\mathrm{nu}}
\renewcommand{\int}{\mathrm{int}}
\newcommand{\ap}{\mathrm{ap}}
\newcommand{\dap}{\mathrm{dap}}
\newcommand{\pre}{\mathrm{pre}}
\newcommand{\conn}{\mathrm{conn}}
\newcommand{\const}{\mathrm{const}}
\newcommand{\conv}{\mathrm{conv}}
\newcommand{\red}{\mathrm{red}}
\newcommand{\falg}{\wedge\mathrm{alg}}

\newcommand{\Spec}{\mathrm{Spec}}
\newcommand{\Alg}{\mathcal{A}lg}
\newcommand{\strAlg}{\mathrm{Alg}}
\newcommand{\CAlg}{\mathcal{CA}lg}
\newcommand{\strCAlg}{\mathrm{CAlg}}
\newcommand{\Mod}{\mathcal{M}od}
\newcommand{\strMod}{\mathrm{Ch}}
\newcommand{\Map}{\mathrm{Map}}
\newcommand{\Art}{\mathrm{Art}}

\newcommand{\Moduli}{\mathrm{Moduli}}

\newcommand{\Ind}{\mathrm{Ind}}

\newcommand{\cl}{\mathrm{cl}}
\newcommand{\gr}{\mathrm{gr}}
\newcommand{\op}{\mathrm{op}}
\newcommand{\fil}{\mathrm{fil}}
\newcommand{\cfil}{\widehat{\mathrm{fil}}}
\newcommand{\calP}{\mathcal{P}}
\newcommand{\calC}{\mathcal{C}}

\newcommand{\calD}{\mathcal{D}}
\newcommand{\calF}{\mathcal{F}}
\newcommand{\calL}{\mathcal{L}}
\newcommand{\calA}{\mathcal{A}}
\newcommand{\calB}{\mathcal{B}}
\newcommand{\calE}{\mathcal{E}}
\newcommand{\calT}{\mathcal{T}}
\newcommand{\rmM}{\mathrm{M}}
\newcommand{\calO}{\mathcal{O}}
\newcommand{\fg}{\mathfrak{g}}
\newcommand{\Br}{\mathrm{Br}}
\newcommand{\SC}{\mathrm{SC}}
\newcommand{\colim}{\operatorname{colim}\limits}
\newcommand{\eps}{\varepsilon}
\newcommand{\ddr}{\mathrm{d_{dR}}}
\newcommand{\eqrightarrow}{\overset{\sim}{\longrightarrow}}
\renewcommand{\cot}{\mathbb{L}}
\newcommand{\sqz}{\mathrm{sqz}}
\renewcommand{\tan}{\mathbb{T}}
\newcommand{\vari}{\textendash}

\newcommand{\defeq}{\mathbin{:\!=}}
\newcommand{\defterm}[1]{\textbf{\emph{#1}}}

\begin{document}
\title{Shifted Lagrangian thickenings of shifted Poisson derived schemes}
\author{Nikola Tomi\'c}
\email{nikola.tomic@umontpellier.fr}
\address{IMAG – UMR 5149 - Université de Montpellier - Place Eugène Bataillon 34090 Montpellier France}
\begin{abstract}
	We prove that the space of shifted Poisson structures on a derived scheme $X$ locally of finite presentation is equivalent to the space of shifted Lagrangian thickenings out $X$, solving a conjecture in shifted Poisson geometry. As a corollary, we show that for $M$ a compact oriented $d$-dimensional manifold and an $n$-shifted Poisson structure on $X$, the mapping stack $\Map(M,X)$ has an $(n-d)$-shifted Poisson structure. It extends a known theorem for shifted symplectic structures to shifted Poisson structures.
\end{abstract}
\maketitle
\tableofcontents
\section{Introduction}
\subsection{Background}
\par Poisson geometry is the study of Poisson manifolds, that is, manifold $M$ with a \defterm{Poisson bracket} on $M$. It is a Lie bracket
\[\{\vari,\vari\} : C^\infty(M)\otimes C^\infty(M) \rightarrow C^\infty(M)\]
which satisfies the Leibniz rule in each variable. The study of such bracket goes back to the work of Lie, Cartan and Dirac. On manifolds, they have been studied by Lichnerowicz (\cite{poisman}). A lot of examples of Poisson manifolds come from physical systems and in fact generalize symplectic manifolds. It is well known since the work of Weinstein (\cite{sympfolman}) that Poisson manifolds locally split into a union of symplectic manifolds. More globally, the Poisson bracket $\{\vari,\vari\}$ gives rise to a vector field $\{f,\vari\}$ for each $f \in C^{\infty}(M)$, called the \defterm{Hamiltonian vector field} defined by $f$. The operation $f \mapsto \{f,\vari\}$ defines a morphism of vector bundles
\[\Theta_{\pi} : T^*M \rightarrow TM.\]
The image of such a morphism is involutive, defining a foliation on $M$. The leaves of this foliation are by nature symplectic, and it gives locally the splitting of the manifold into union of symplectic manifolds. Reciprocally, the data of a symplectic foliation on a manifold provides a Poisson bracket on it (see, for instance, the book of Vaisman (\cite{lectpoisman}) for an exposition of this topic). It is possible to describe the Poisson structure in another way by integrating the foliation into a Lie groupoid over $M$. We get an object called a symplectic groupoid, and such data recovers the Poisson bracket (see, for instance, the work of Weinstein (\cite{sympgrpdman})).
\par We focus here on an analog of Poisson geometry within the setting of Derived geometry. Derived geometry is a generalization of algebraic geometry and allows people to deal with singularities in a convenient way. The theory was first given by To\"en and Vezzosi (\cite{HAGII}), Lurie (\cite{DAG}, but also \cite{SAG}, which is still not fully written yet), and Gaitsgory and Rozenblyum (\cite{GRI,GRII}, some unproven claims in the first book were recently proven by different authors, see \cite{endGR} for a review). Since then, there have been a lot of applications of the theory. One outstanding application is the theory of shifted symplectic structures, as defined by Pantev, To\"en, Vaqui\'e and Vezzosi (\cite{PTVV}), providing the existence of many new structures on moduli stacks and allowing the construction of new invariants. For instance, the categorification of Donaldson-Thomas invariants by Brav, Bussi, Dupont, Joyce, and Szendr\H oi (\cite{BBDJS}) has been done using ideas from shifted symplectic geometry. Recently, the construction has been generalized by Hennion, Holstein, and Robalo (\cite{HHRI,HHRII}).
\par Shifted Poisson geometry is the derived analog of Poisson geometry. Calaque, Pantev, To\"en, Vaqui\'e, and Vezzosi (\cite{CPTVV}), and, independently, Pridham (\cite{poisndpri}), define shifted Poisson structures on derived Artin stacks and prove that non-degenerate shifted Poisson structures are equivalent to shifted symplectic structures. This last fact is easy to prove in differential geometry, but it is non-trivial in derived geometry, specifically because one has to provide a homotopy Poisson algebra structure given a symplectic structure, which involves an infinite number of operations. Another thing that is hard to establish is the correct notion of Poisson structures for Artin stacks. The naive approach does not work because shifted polyvector fields, used to define shifted Poisson structures, are not functorial, so we can't locally define shifted Poisson structures, it is however functorial on formally etale maps so it is fine for Deligne--Mumford stacks. In order to solve this problem, the authors of \cite{CPTVV} develop the technique of formal localization. What is really hidden in their paper is the use of deformation theory of operads. The deformation theory for operads has been developed by Ginzburg and Kapranov (\cite{GKKoz}), and Loday and Vallette (\cite{LV}). The theory is used extensively by Melani (\cite{poisbiv}), and later by Melani and Safronov (\cite{MSI}) in order to define coisotropic structures, the ``Lagrangian" counterpart for morphisms into a shifted Poisson derived stack. In a converse point of view, shifted symplectic groupoids have been studied recently by Calaque and Safronov (\cite{sympgrpd}). They characterize shifted symplectic groupoids on $X$ in terms of Lagrangian morphisms out of $X$.
\subsection{Main results}
\par We have seen that each Poisson manifold has a symplectic foliation, and the bracket is totally determined by the symplectic foliation. In the derived setting, a similar phenomenon appears. First, one has to adapt the notion of foliation in the derived case by looking at \defterm{derived foliations}, firstly investigated by To\"en and Vezzosi (\cite{fol}). Because we work algebraically, it is not possible to integrate globally a derived foliation anymore, we are able to only formally integrate around a point. Let $X$ be an $n$-shifted Poisson derived scheme, thanks to the results of \cite{CPTVV}, such a structure can be described as a weight greater than $2$ Maurer--Cartan element
\[\pi \in \Pol(X,n) \defeq \widehat\Sym(\tan_X[-n]).\]
It is an element satisfying the equation
\[\frac12[\pi,\pi] + \d\pi = 0\]
where the Lie bracket is the Schouten--Nijenhuis bracket. In that case, twisting the differential of $\Pol(X,n)$ by $[\pi,\vari]$, and remembering the symmetric grading will give a mixed graded algebra over $X$ that is actually a derived foliation. We denote it by $\calF$. If we follow geometric intuition, the induced map
\[\Theta_\pi \colon \cot_X[n] \rightarrow \tan_X.\]
is the map from the tangent complex of the foliation to the one of $X$. We cannot expect this foliation to be symplectic since this would imply that the cotangent complex of $X$ is equivalent to the tangent complex up to a shift, which is not the case in general. However, we can look at the tangent complex of the leaves derived prestack, which will be a formal derived prestack in our case. Thanks to the results of Brantner, Magidson, and Nuiten (\cite{deffol}) (some of the results in this paper were previously studied by Nuiten in a more restrictive context, and using tools of model categories (\cite{KozLie}), and by Calaque and Grivaux (\cite{surveyform})), and Fu (\cite{follie}), when $X$ is a locally coherent derived scheme, we know that there is a formal thickening,
\[\pi : X \rightarrow \left[X/\calF\right],\]
and a fiber sequence:
\[\cot_X[n] \rightarrow \tan_X \rightarrow \pi^* \tan_{\left[X/\calF\right]}.\]
If we dualize this sequence, we get the sequence
\[\pi^*\cot_{\left[X/\calF\right]} \rightarrow \cot_X \rightarrow \tan_X[-n],\]
thus giving an equivalence
\[\pi^*\cot_{\left[X/\calF\right]}[n+1]\simeq \pi^*\tan_{\left[X/\calF\right]}.\]
So, one may expect $\left[X/\calF\right]$ to be $(n+1)$-shifted symplectic with a Lagrangian structure on $\pi$ (since $\tan_{\pi} \simeq \cot_X[n]$). In this way, we have constructed what is called a \defterm{$(n+1)$-shifted Lagrangian thickening} $X \rightarrow \calL$.
\par Reciprocally, after the results of Melani and Safronov (\cite{MSI,MSII}), we might expect that a Lagrangian structure is a special case of what is called a coisotropic structure, and such a structure gives a Poisson structure on the source. That is how we should recover the original Poisson structure. We state the main result of this paper:
\begin{theorem}(\ref{thm:main})
	Let $X$ be a locally of finite presentation derived scheme, denote by $\Pois(X,n)$ the space of $n$-shifted Poisson structures on $X$, and by $\LagThick(X,n+1)$ the space of $(n+1)$-shifted Lagrangian thickenings. There is an equivalence:
	\[\Pois(X,n) \simeq \LagThick(X,n+1)\]
\end{theorem}
From this theorem, we can deduce the existence of new Lagrangian and Poisson structures.
\begin{corollary}(\ref{cor:aksz})
	Let $Y$ be a $\calO$-compact with $\calO$-orientation of degree $d$ derived stack, in the sense of \cite[2.1.]{PTVV}. Let $X$ be a locally of finite presentation derived scheme with an $n$-shifted Poisson structure. Then there exists a map $\Map(Y,X) \rightarrow \calL$ of derived prestacks with an $(n-d+1)$-shifted Lagrangian structure. If $Y = M_{B}$ is the Betti stack of a compact oriented $d$-dimensional manifold, then there exists an $(n-d)$-shifted Poisson structure on $\Map(M_B,X)$.
\end{corollary}
Given a Poisson structure, the existence of a Lagrangian thickening has already been conjectured by Safronov in his lecture notes (\cite[3.2.]{lectpois}). In these notes, formal thickening is called the ``symplectic realization." For the converse direction, given a Lagrangian thickening $X \rightarrow Y$ between two derived prestacks, one may want to use the results of Melani and Safronov about coisotropic structures. However, they prove that Lagrangian morphisms between Artin stacks have a coisotropic structure and then a Poisson structure on the source, which is not the right hypothesis in our setting, $X$ is a derived scheme, but $Y$ may fail to be geometric (it is currently unknown whether or not a formal thickening of a derived scheme satisfies etale descent). We bypass this by proving the result before taking the formal thickening. We summarize the steps of the proof of the theorem as follows:
\begin{itemize}
	\item A first step is to provide the Lagrangian structure on $\pi$. It comes from the fact that whenever the map $\calF\rightarrow \calO_X$ has an internal shifted Lagrangian structure of weight $-1$, then $[X/\calF]$ automatically has a shifted symplectic structure, and $X \rightarrow [X/\calF]$ has a shifted Lagrangian structure. Reciprocally, if $X \rightarrow [X/\calF]$ has a shifted Lagrangian structure, then $\calF \rightarrow \calO_X$ as an internal shifted Lagrangian structure of weight $-1$ (\cref{thm:lagalglagthick}). This will be investigated in \cref{sec:firststep}.
	\item Now, it suffices to provide a shifted Lagrangian structure of weight $-1$ on the foliation $(\Pol(A,n),\{\pi,\vari\})$. This is exactly given by the Schouten--Nijenhuis bracket, which is non-degenerate. It then gives a shifted symplectic structure of weight $-1$ (the bracket is of weight $-1$). In order to properly state that, we have to regive all the deformation theoretic features developed by Melani and Safronov in \cite{MSI} in the mixed graded setting, so we have to study the operad $\bP_n^{\gr}$ of Poisson algebras with bracket of degree $1-n$ and of weight $-1$. The model category of mixed graded complexes is not well-behaved to do deformation theory (because either the unit is not cofibrant, if we choose the projective model structure, either one has to consider fibrant resolutions, if we choose the injective model structure), and so we model the underlying $\infty$-category of graded mixed complexes by using complete filtered complexes instead. The fact that the category of complete filtered complexes models the same $\infty$-category as mixed graded complexes is a known result that comes from Koszul duality between $k[\eps]$ and $k[\hbar]$. It is proven in a model-independent setting for any stable category after Ariotta (\cite{filcohchain}). The deformation theory of operads in that setting has been investigated by Calaque, Campos, and Nuiten (\cite{curvlie}), and it is well-behaved. This step will be investigated in \cref{sec:sndstep}.
	\item The last thing to do is to establish the affine case and descent for Lagrangian thickenings and Poisson structures. We firstly observe that the space of shifted Poisson structures on $\calO_X$ is equivalent to the space of shifted Poisson structures of weight $-1$ on $\calO_X(0)$ (mixed graded $\calO_X$ concentrated in weight $0$). Then, given a shifted Lagrangian structure of weight $-1$ on a morphism targeting $\calO_X(0)$, we know by an analog result of \cite{MSI} that there is a shifted Poisson structure of weight $-1$ on $\calO_X(0)$ and then a shifted Poisson structure on $X$. That way we bypass the difficulty that we had because of the non-geometricity of our derived prestacks. Zariski descent for Poisson structures is a known result (more generally it satisfies etale descent). It comes from the fact that the tangent complex is functorial along formally etale morphisms. However, Zariski descent for Lagrangian thickenings is not known, but it is an almost direct consequence of Zariski descent for formal thickenings, which was proven in \cite{deffol}. Finally, we prove the main theorem and its corollary. We outline those steps in \cref{sec:thrdstep}.
\end{itemize}
We mention that Melani and Pourcelot are currently investigating a more explicit and general approach to show the existence results of \cref{cor:aksz}, and that Johnson-Freyd has previously established a similar result in \cite{akszjohnfreyd}.
\section{Preliminaries}
\par In this paper, we assume that the word ``category" refers to an $\infty$-category unless specified otherwise.\\
\par We fix $k$, a field of zero characteristic. In that case we fix the following notations:
\begin{itemize}
	\item $\Mod_k$ is the $\infty$-category of cochain complexes of $k$-modules,
	\item $\CAlg_k$ is the $\infty$-category of commutative differential graded $k$-algebras,
	\item $\CAlg_k^{\conn}$ is the $\infty$-category of connective commutative differential graded $k$-algebras.
\end{itemize}
\par Recall that $\Mod_k$ has a t-structure; we may write $\Mod_k^{\conn}$ for the category of connective modules. Because the content of \cite{deffol} is settled up in the context of the theory of derived rings of \cite{derdr}, we precise that ``animated" in their setting means ``connective" and ``derived" is irrelevant because we work in characteristic zero. In the sequel we will also denote by $\calC$ any $k$-linear presentably symmetric monoidal category (by $k$-linear we mean stable tensored over $\Mod_k$), the tensor product will be usually denoted $\otimes$. If $A$ is an algebra in $\calC$, the relative tensor product will be denoted $\vari \otimes_A \vari $ and the $A$-linear internal Hom will be denoted $\underline\Hom_A(\vari,\vari)$.
\subsection{Graded, filtered, and complete filtered objects}
\par We recall here basic facts about filtered and graded objects with values in a stable category, we refer to \cite[2.6. \& 3.1.]{luriefil} or \cite{modfil} for proofs of the claimed statements.
\par Let $\Z$ be the poset of integers with the order $\leq$. The category $\calC^{\fil(, \geq 0)}$ of (non-negatively) \defterm{filtered objects} is defined as follows:
\[\calC^{\fil(, \geq 0)} \defeq \Fun\left((\Z_{(\geq 0)},\leq)^{\op},\calC\right).\]
An object $\calC^{\fil}$ will be denoted $F^{\bullet} M$ with maps
\[\dots \rightarrow F^i M \rightarrow F^{i-1}M\rightarrow \dots.\]
An object of $\calC^{\fil, \geq 0}$ will be denoted the same way, but the filtration is bounded on the right:
\[\dots \rightarrow F^1M \rightarrow F^0M.\]
We will denote by $M(i)$ the filtered object given by
\[F^j (M(i)) \defeq \begin{cases} 0 & \text{ if } j > i\text{,}\\ M & \text{elsewhere.}\end{cases}.\]
We set
\[\calC^{\gr(,\geq 0)} \defeq \Fun\left(\Z^{\delta}_{(\geq 0)},\calC\right)\]
for the category of (non-negatively) \defterm{graded objects} of $\calC$, where $\Z^{\delta}$ is the discrete integer category. An object of $\calC^{\gr(, \geq 0)}$ will be denoted $M^\bullet$. There is a functor denoted $\gr$:
\[\begin{array}{ccc}
	\calC^{\fil(,\geq 0)} & \longrightarrow & \calC^{\gr(, \geq 0)}\\
	F^\bullet M & \longmapsto & \left(\cofib(F^{i+1} M \rightarrow F^i M)\right).
\end{array}\]
that is called the associated graded functor. If $M \in \calC$ and $i \in \Z$ we will denote by $M(i)$ the graded object of $\calC$ $M$ put in weight $i$ and $0$ elsewhere. The notation is the same as the one in the filtered case; however, it will be clear from the context what object we are talking about. We can see that it is also compatible with the functor $\gr$. $\calC^{\fil(, \geq 0)}$ and $\calC^{\gr(, \geq 0)}$ are $k$-linear and presentable, and $\gr$ commutes with limits and colimits. Using the addition of $\Z$ and $\Z^\delta$, we get presentably symmetric monoidal products on both of them by Day convolution (\cite[2.2.6.]{HA}).
\par Recall that an object $F^{\bullet} M \in \calC^{\fil}$ is called \defterm{complete} if $\lim_n F^n M \simeq 0$. We denote $\calC^{\cfil}$ and $\calC^{\cfil,\geq 0}$ the full subcategories spanned by complete filtered objects. Remark that $M(i)$ is always complete by definition. There is an adjunction:
\[\widehat{(\vari)} : \calC^{\fil(, \geq 0)} \rightleftarrows \calC^{\cfil(, \geq 0)} : i\]
where $i$ is the inclusion and $\widehat{(\vari)}$ is the completion. It is actually a left exact localization at morphisms inducing equivalences on associated graded. $\calC^{\cfil(, \geq 0)}$ is then $k$-linear presentably symmetric monoidal (The symmetric monoidal product being given by the completion of the monoidal product of filtered objects). The functor $\gr$ factors as
\[\gr : \calC^{\cfil(, \geq 0)} \rightarrow \calC^{\gr(,\geq 0)}.\]
It is conservative, symmetric monoidal, and has right and left adjoints. It induces functors
\[\CAlg(\calC^{^{^{(\!}}\cfil^)(,\geq 0)}) \rightarrow \CAlg(\calC^{\gr(, \geq 0)})\]
that also have right and left adjoints (and is conservative in the complete case). It is also well-behaved with the internal Hom functors in the complete case.
\begin{lm}\label{lm:grhom}
	Let $A \in \CAlg(\calC^{\cfil})$ and $X,Y \in \Mod_A(\calC^{\cfil})$. There is a natural equivalence
	\[\gr\underline\Hom_A(X,Y) \simeq \underline\Hom_{\gr A}(\gr X,\gr Y).\]
\end{lm}
\begin{proof}
	Let $Z \in \Mod_{\gr A}(\calC^{\gr})$. We have
	\begin{align*}
		\Map_{\gr A}(Z,\gr\underline\Hom_A(X,Y)) & \simeq \Map_A(LZ,\underline\Hom_A(X,Y))\\
							 & \simeq \Map_A(LZ\otimes_A X,Y)\\
							 & \simeq \Map_A(L(Z\otimes_{\gr A} \gr X), Y)\\
							 & \simeq \Map_{\gr A}(Z\otimes_{\gr A} \gr X, \gr Y)\\
							 & \simeq \Map_{\gr A}(Z, \underline\Hom_{\gr A}(\gr X, \gr Y)).
	\end{align*}
	Where $L$ is the left adjoint of $\gr$. It is given by
	\[(Z^i) \mapsto \dots \overset{0}{\rightarrow} F^j(LZ) \defeq Z^{j-1}[-1] \overset0\rightarrow \dots,\]
	with an $A$-module structure coming from the quotient $A \mapsto \gr A$. The equivalence
	\[LZ \otimes_A Y \simeq L(Z\otimes_{\gr A} \gr X)\]
	is proved as follows: If we forget about the $A$-linear structure, we have
	\[\left(\gr \left(L Z\otimes X\right)\right)^n \simeq \bigoplus_{i+j = n}(Z^{i-1}[-1]\oplus Z^i)\otimes (\gr X)^j\simeq \left(\gr L(Z\otimes\gr X)\right)^n\]
	giving
	\[L Z \otimes X \simeq L(Z\otimes \gr X).\]
	If we take into account the $A$-linear structure, we also have
	\[L Z \otimes_A X \simeq L(Z \otimes_{\gr A} X)\]
	because $L$ also commutes with colimits. By Yoneda,
	\[\gr\underline\Hom_A(X,Y) \simeq \underline\Hom_{\gr A}(\gr X, \gr Y).\]
\end{proof}
\par If $G : \calC \rightarrow \calD$ is lax monoidal, where $\calD$ satisfies the same assumptions as $\calC$, then we have induced lax monoidal functors
\[G^{\fil(, \geq 0)} : \calC^{\fil(, \geq 0)} \rightarrow \calD^{\fil(, \geq 0)}\]
and
\[G^{\gr(, \geq 0)} : \calC^{\gr( \geq 0)} \rightarrow \calD^{\gr(, \geq 0)}.\]
If it is moreover symmetric monoidal and commutes with colimits, then the induced lax structure is symmetric monoidal because the Day convolution tensor product is a left Kan extension. If $G$ commutes with limits, there is an induced lax monoidal functor
\[G^{\cfil(, \geq 0)} : \calC^{\cfil(, \geq 0)} \rightarrow \calD^{\cfil(, \geq 0)}.\]
There is an ``underlying object" symmetric monoidal colimit-preserving functor
\[(\vari)^u : \calC^{\fil(, \geq 0)} \rightarrow \calC\]
sending $F^\bullet M$ to $\colim_n F^n M$. It also induces a left adjoint on the level of algebras:
\[(\vari)^u : \CAlg(\calC^{\fil(, \geq 0)}) \rightarrow \CAlg(\calC).\]
Note that in the non-negatively filtered case, these functors send $F^\bullet M$ to $F^0 M$. In the complete case it is slightly different. The underlying object functor is neither a right adjoint nor a left adjoint. However, the functor $F^0(\vari)$ is lax monoidal and has a symmetric monoidal left adjoint given by $(\vari)(0)$. By tradition (see \cite[1]{CPTVV}, where complete filtered objects are seen as mixed graded objects instead in shifted symplectic geometry), we will denote $F^0(\vari)$ by $|\vari|$ for complete filtered objects. The functor $(\vari)(0)$ also has a left adjoint traditionally called ``left realization". Explicitly, it sends $F^\bullet M$ to $\cofib\left(F^1 M \rightarrow (F^\bullet M)^u\right)$.
\par When $\calC$ has a t-structure, $\calC^{\fil(, \geq 0)}$, $\calC^{\gr(, \geq 0)}$, and $\calC^{\cfil(, \geq 0)}$ have t-structures as well called the neutral one defined, for instance, in \cite[3.3.5.]{derdr}. For instance, for a filtered object $X$, $X$ is connective in this t-structure if and only if $(\gr X)^i$ is connective for each $i$. If the initial t-structure was right complete, the neutral t-structures are right complete as well.
\subsection{Adic and de Rham filtrations}\,
\par There is a functor 
\[\aug : \CAlg(\calC^{\fil, \geq 0}) \rightarrow \CAlg(\calC)^{\Delta^1}\]
sending $F^\bullet M$ to the map $F^0 M \rightarrow \gr(F^\bullet M)^0$. It is accessible, commutes with limits, and then has a left adjoint called $\adic$ that we call the \defterm{adic filtration functor}. In the complete case, the same functor
\[\aug : \CAlg(\calC^{\cfil, \geq 0}) \rightarrow \CAlg(\calC)^{\Delta^1}\]
exists and also has a left adjoint that is the completion of $\adic$. It is denoted on $B \rightarrow A$ by $\DR(A/B)$. It is called the \defterm{relative de Rham algebra} of $B \rightarrow A$. Recall from \cite[7.3.]{HA} that there is a cocartesian fibration $T_{\CAlg(\calC)} \rightarrow \CAlg(\calC)$ classifying the covariant functor $A \mapsto \Mod_A$. Hence we have the cotangent complex functor defined as a left adjoint to the square zero extension functor:
\[\cot_{(\vari)} : \CAlg(\calC) \rightleftarrows T_{\CAlg(\calC)} : \sqz.\]
It applies for $\calC$, but also for $\calC^\fil$, $\calC^\gr$, $\calC^{\cfil}$... Note that because we are in zero characteristic, it coincides with the cotangent complex formalism from \cite{derdr}. We can relate the adic and de Rham filtrations to the cotangent complex:
\begin{lm}\label{lm:grdr}
	Let $B\rightarrow A$ be an object of $\CAlg(\calC)^{\Delta^1}$, then we have the following identities:
	\begin{itemize}
		\item $\gr(\adic(B\rightarrow A)) \simeq \Sym_A\left(\cot_{A/B}(1)[-1]\right)$ and
		\item $\gr(\DR(A/B)) \simeq \Sym_A\left(\cot_{A/B}(1)[-1]\right)$.
	\end{itemize}
\end{lm}
\begin{proof}
	It is because the following diagram of right adjoints commutes:
	\[\begin{tikzcd}
		\CAlg(\calC)^{\Delta^1} & \CAlg(\calC^{^{^(\!}\cfil^), \geq 0}) \ar[l,"\aug"]\\
		\Mod(\calC) \ar[u,"\delta\sqz"] & \CAlg(\calC^{\gr, \geq 0}) \ar[u,"R"] \ar[l,"\ev_{01}"]
	\end{tikzcd}\]
	where:
	\begin{itemize}
		\item $\Mod(\calC)$ is the category of pairs $(A,M)$ with $A \in \CAlg(\calC)$ and $M \in \Mod_A(\calC)$,
		\item $\ev_{01}$ sends $A^\bullet$ to $(A^0,A^1[1])$. Its left adjoint is $(A,M) \mapsto \Sym_A(M(1)[-1])$,
		\item $\delta\sqz$ sends $(A,M)$ to $A\rightarrow\sqz(A,M)$. Its left adjoint is the relative cotangent complex $\cot_{\vari/\vari}$, and
		\item $R$ sends $A^\bullet$ to the filtered algebra $\dots \overset{0}\rightarrow A^1 \overset{0}\rightarrow A^0$. Its left adjoint is $\gr$ because a map $F^\bullet X \rightarrow R(A^\bullet)$ is the same data as the data of maps $F^i X \rightarrow A(i)$ with a homotopy from $F^{i+1} X \rightarrow F^i X \rightarrow A(i)$ to zero. That is, exactly a map $(\gr X)^i \rightarrow A(i)$.
	\end{itemize}
	The diagram of left adjoints commutes. 
\end{proof}
The adic filtration of $B \rightarrow A$ is really a filtration on $B$ in the following sense:
\begin{lm}\label{lm:adicfil}
	For $B \rightarrow A$, an object of $\CAlg(\calC)^{\Delta^1}$, we have
	\[F^0\left(\adic(B\rightarrow A)\right) \simeq B.\]
\end{lm}
\begin{proof}
	Thanks to \cref{lm:grdr}, we have that
	\[F^0(\adic(B\rightarrow A) \simeq (\adic(B\rightarrow A))^u.\]
	But the right adjoint of $(\vari)^u$ is the constant diagram functor $C \mapsto (\dots \rightarrow C \rightarrow C)$. Composition with augmentation gives $C \mapsto (C \rightarrow 0)$, which is right adjoint to $(B \rightarrow A) \mapsto B$.
\end{proof}
\subsection{The de Rham functor}\,
\par We have defined the relative de Rham functor. It is possible to define the absolute de Rham functor as a special case of the relative one, for $A\in\CAlg(\calC)$, we set
\[\DR(A) \defeq \DR(A/\mathbf{1})\]
where $\mathbf{1}$ is the unit of $\calC$. It is left adjoint to the $\gr^0$ functor $\CAlg(\calC^{\cfil})\rightarrow\CAlg(\calC)$. Of course, we still have
\[\gr\DR(A) \simeq \Sym\left(\cot_A(1)[-1]\right).\]
We shall compare different de Rham functors that are computed in different categories. We set $\calC'\defeq\calC^{\cfil}$ and $\calC''\defeq\calC^{\cfil,\geq 0}$. For convenience, we will write $\DR^{\int}$ for the de Rham functor computed for algebras object in $\calC'$ or $\calC''$. We will need the following properties:
\begin{prop}\label{prop:drcomtrivgr}
	For $A \in \CAlg(\calC)$, $\DR(A) \in \CAlg(\calC^{\cfil})$, and the following diagram:
	\[\begin{tikzcd}
		\CAlg(\calC) \ar[r,"(\vari)(0)"] \ar[d,"\DR"] & \CAlg(\calC')\ar[d,"\DR^{\int}"]\\
		\CAlg(\calC^{\cfil}) \ar[r,"{(\vari)(0)}^{\int}"] & \CAlg(\calC'^{\cfil})
	\end{tikzcd}\]
	commutes, where ${(\vari)(0)}^{\int}$ is obtained from $(\vari)(0) : \calC \rightarrow \calC'$ by taking ${(\vari)(0)}^{\cfil} : \calC^{\cfil} \rightarrow \calC'^{\cfil}$ and inducing it on commutative algebras since this functor is lax monoidal and preserves limits.\\
	For $B \in \CAlg(\calC'')$, the following diagram:
	\[\begin{tikzcd}
		\CAlg(\calC'') \ar[r,"\gr^0"] \ar[d,"\DR^{\int}"] & \CAlg(\calC)\ar[d,"\DR"]\\
		\CAlg(\calC''^{\cfil}) \ar[r,"\gr^{0,\int}"] & \CAlg(\calC^{\cfil})
	\end{tikzcd}\]
	commutes. Here $\gr^{0,\int}$ is defined in a similar way as ${(\vari)(0)}^{\int}$.
\end{prop}
\begin{proof}
	It's enough to check the commutativity on right adjoints:
	\[\begin{tikzcd}
		\CAlg(\calC) & \CAlg(\calC')\ar[l,"|\vari|"]\\
		\CAlg(\calC^{\cfil}) \ar[u,"\gr^0"] & \CAlg(\calC'^{\cfil}), \ar[u,"\gr^0"] \ar[l,"|\vari|^{\int}"]
	\end{tikzcd}\]
	and
	\[\begin{tikzcd}
		\CAlg(\calC'') & \CAlg(\calC)\ar[l,"(\vari)(0)"]\\
		\CAlg(\calC''^{\cfil}) \ar[u,"\gr^0"] & \CAlg(\calC^{\cfil}).\ar[u,"\gr^0"] \ar[l,"(\vari)(0)^{\int}"]
	\end{tikzcd}\]
	They commute because they come from commutative diagrams on the underlying $k$-linear categories.
\end{proof}
The remaining subsections are irrelevant for \cref{sec:firststep}.
\subsection{Model categorical settings}\label{sec:modelcat}\,
\par We will always write $\rmM$ for a model category satisfying the standing assumptions of \cite[1.1.]{CPTVV}. That is, it is $k$-linear and has a monoidal model category structure. For instance, the model category of cochain complexes $k$-vector spaces $\strMod_k$ satisfies the standing assumptions, and we have $\strMod_k[Qiso^{-1}]\simeq \Mod_k$. The assumptions are made in such a way that the localization $\calC \defeq \rmM[W_{eq}^{-1}]$ is $k$-linear presentably symmetric monoidal. When $A$ is a algebra object of $\rmM$, we write $\strMod_A(\rmM)$ for the model category of cochain complexes of $A$-modules. It localizes into the $\infty$-category $\Mod_A(\calC)$. It also satisfies the standing assumptions. In \cite{MSI}, the authors are usually interested in the realization functor $|\vari| : \calC \rightarrow \Mod_k$ coming from the cotensor structure over $\Mod_k$. But here, we will be interested in the filtered analog of it.
\par Denote by $\rmM^{\cfil}$ the $1$-category of complete filtered objects, that is, sequences of injections:
\[\cdots \hookrightarrow F^{1}C \hookrightarrow F^0 C \hookrightarrow F^{-1}C \hookrightarrow \cdots,\]
where $C \defeq \colim F^i C$ and $\lim C/F^i C \simeq C$. There is a functor $\rmM^{\cfil} \rightarrow \rmM^{\gr}$ given by $F^{\bullet} C \mapsto \left(F^i C/ F^{i+1}C\right)$. 
\par The completed tensor product of two complete filtered objects makes $\rmM^{\cfil}$ into a symmetric monoidal model category (\cite[2.20]{curvlie}). By \cite[2.19.]{curvlie}, it satisfies the standing assumptions. By \cite[4.3.]{modfil}, $\rmM^{\cfil}[W_{eq}^{-1}] \simeq \calC^{\cfil}$ so we really get a model for complete filtered objects of $\calC$. $\calC^{\cfil}$ is naturally tensored and cotensored over $\Mod_k^{\cfil}$ and we will consider the filtered realization functor:
\[|\vari|^{\cfil} : \calC^{\cfil} \rightarrow \Mod_k^{\cfil}.\]
This functor can be explicitly computed as the derived functor $\mathbb{R}\Hom(\mathbf{1}_{\rmM^{\cfil}},\vari)$ (one only have to consider a fibrant replacement in the second variable in order to compute this since the unit is supposed to be cofibrant). Recall from \cite[Section 1.2.]{MSI} mixed graded structures:
\begin{defn}
	A \textbf{mixed graded object of $\rmM$} is the data of objects $X^i \in \rmM$ for $i \in \Z$ together with morphisms $\eps : X^{\bullet} \rightarrow X^{\bullet+1}[1]$ satisfying $\eps^2 = 0$. We denote by $\rmM^{\eps,\gr}$ the category of mixed graded objects of $\rmM$. $\rmM^{\eps,\gr}$ can be endowed with the injective model structure. In that case, it satisfies the standing assumptions. We will denote by $\calC^{\eps,\gr}$ the associated $\infty$-category.
\end{defn}
\par For complexes, there is a weak version of it that can be found by taking a quasi-free resolution of $k[\eps]$. See \cite[2.2.1.]{curvlie} for a discussion on it.
\begin{defn}
	A \textbf{weak mixed graded complex} is the data of $k$-linear chain complexes $X^i$ for $i \in \Z$ together with morphisms $\eps_k : X^{\bullet} \rightarrow X^{\bullet+k}[1]$ satisfying $\sum_{i+j=n} \eps_i \circ \eps_j = 0$ ($\eps_0 = d$). The category of weak mixed graded complexes gives the same $\infty$-category $\Mod_k^{\eps,\gr}$ as before. This definition will be relevant because some objects appearing later are naturally weakly mixed graded instead of strictly.
\end{defn}
\par In general, we have that the inclusion functor
\[\calC \rightarrow \calC^{\eps,\gr}\]
has both right and left adjoints. We denote by $|\vari|^r$ the right adjoint and $|\vari|^l$ the left adjoint. In practice, when we consider twisting by a differential, we mean that we can endow our twisted object by a mixed graded structure, and then take the right realization or the left realization depending whether the mixed structure is completely in non-negative or non-positive degrees.
\begin{remark}
	$\calC^{\eps,\gr}$ is actually equivalent to $\calC^{\cfil}$. A model-independent proof is known for any stable category in \cite{filcohchain}. In a more explicit way, this is the incarnation of Koszul duality between $k[\eps]$ and $k[\hbar]$. Indeed, we may identify $\rmM^{\eps,\gr}$ with the model category $\strMod_{k[\eps]}(\rmM^{\gr})$ where $\eps$ has weight $1$ and cohomological degree $1$. Koszul duality gives a functor
	\[\strMod_{k[\eps]}(\rmM^{\gr}) \rightarrow \strMod_{k[\hbar]}(\rmM^\gr)\]
	where $\hbar$ has weight $-1$ and cohomological degree $0$. This functor is then fully faithful on the associated $\infty$-categories with essential image objects $M^i$ such that
	\[\lim \left(\dots \overset\hbar\rightarrow M^i \overset\hbar\rightarrow M^{i-1} \overset\hbar\rightarrow \dots \simeq 0\right),\]
	which is clearly $\rmM[W_{eq}^{-1}]^{\cfil}$. A proof of this for complexes over a field has been done in \cite[2.27.]{curvlie}. If we write stuff explicitly, the equivalence $\calC^{\eps,\gr} \simeq \calC^{\cfil}$ is given by
	\[X \mapsto (\dots \rightarrow |X\otimes k(-1)|^r \rightarrow |X|^r \rightarrow \dots).\]
	From this perspective, $|\vari|^r$ on $\calC^{\cfil}$ is $F^\bullet M \mapsto F^0 M$ and $|\vari|^l$ is $F^\bullet M \mapsto \cofib(F^1 M \rightarrow (F^\bullet M)^u)$. 
\end{remark}
\subsection{Complete filtered operads and cooperads}\,
\par We define operads and cooperads of complete filtered complexes as the ones that were investigated in \cite{curvlie}.
\begin{defn}
	A \defterm{complete operad} $\calP$ is an operad on $\strMod_k^{\cfil}$. To avoid confusion, for a symmetric sequence $\calP$, we write $\calP((n))$ for the space of arity $n$ operations. Explicitly, it is a symmetric sequence $\calP$ given with operations:
	\[\widehat{\bigoplus_{k\geq 0}} \calP((k)) \hat\otimes_{\Sigma_k} \left(\widehat{\bigoplus_{i_1+\cdots + i_k = n, i_l \geq 0}} \operatorname{Ind}^{\Sigma_n}_{\Sigma_{i_1}\times\cdots\times\Sigma_{\i_k}}\left(\calP((i_1))\hat\otimes\cdots\hat\otimes\calP((i_k))\right)\right) \rightarrow \calP((n))\]
	Here $\widehat\oplus$ and $\widehat\otimes$ are, respectively, the coproduct and the monoidal product of complete filtered complexes. A \defterm{complete cooperad} $\calD$ is a cooperad on $\strMod_k^{\cfil}$.
\end{defn}
For instance, if $A$ is complete filtered complexes, $\End_A$ is naturally filtered and endowed with a complete operad structure. The linear dual $\End_A^{\vee}$ (inside of complete filtered complexes) is a complete cooperad. It is possible to consider complete colored operads in a similar way.
\begin{defn}
	A \defterm{complete colored operad} $\calP$ is a colored operad on the category of complete filtered complexes.
\end{defn}
We will also need the notion of Hopf operads/cooperads; there is the complete Hadamard product of two symmetric sequences $\calP$ and $\calP'$,
\[\left(\calP \hat\otimes_{\mathrm{H}} \calP'\right)((n)) \defeq \calP((n)) \hat\otimes\calP'((n)).\]
One can check that if $\calP$ and $\calP'$ are complete operads, then $\calP \hat\otimes_{\mathrm{H}} \calP'$ is a complete operad endowing the category of complete operads with the structure of a monoidal category, with $\Com$ as the unit.
\begin{defn}
	A \defterm{Hopf complete operad} is a complete operad $\calP$ that is a coalgebra in the monoidal category of complete operads with Hadamard product.
\end{defn}
\par The same thing holds for cooperads. If $\calD$ and $\calD'$ are two complete cooperads, then $\calD \hat\otimes_{\mathrm{H}} \calD'$ is a cooperad.\\
If $\calD$ is a cooperad, we set $\calD^{\cu}$ the symmetric sequence given by $\calD$ in positive arities and $\calD^{\cu}((0)) \defeq \calD((0)) \oplus k$.
\begin{defn}
	A \defterm{unital Hopf complete cooperad} is the data of a complete cooperad $\calD$ such that $\calD^{\cu}$ that is an algebra in the category of complete cooperads with Hadamard product, the unit is the inclusion of the second factor of $\calD((0)) \oplus k$ and $\calD^{\cu}\rightarrow \calD$ is a morphism of complete cooperads.
\end{defn}
In order to handle operad shifts, we introduce a new notation: 
\begin{defn}
	If $\calP$ is a complete filtered symmetric sequence and $\calL$ is a complete invertible complex, we let
	\[\calP\{\calL\} \defeq \calP\hat\otimes_{\mathrm{H}} \End_{\calL}.\]
\end{defn}
This notation coincides with the classical one if we shift by $k[n]$. For instance, $\calP\{n\} = \calP\{k[-n]\}$ and preserves the operadic and cooperad structures. For the operadic structure, it is straightforward, because $\End_{\calL}$ is an operad. For the cooperadic structure, it is because $\End_{\calL} \simeq \End_{\calL^{\vee}}^{\vee}$ which is a cooperad. The same definition for an algebra over an operad holds:
\begin{defn}
	Let $\calP$ be a complete operad. A \defterm{$\calP$-algebra} is a complete filtered complex $A$ with a morphism $\calP \hat\circ A \rightarrow A$, where $A$ is seen as the constant symmetric sequence with $A((i)) = A$ if $i=1$, $0$ otherwise.
\end{defn}
The category of $\calP$-algebras has a structure of model category according to \cite[2.15]{curvlie}
\begin{theorem}
	Let $\calP$ be a complete operad. the category of $\calP$-algebras has a model structure given by
	\begin{itemize}
		\item Weak equivalences are maps inducing a quasi-isomorphism on associated graded.
		\item Fibrations are maps that are surjections in each filtration degree.
	\end{itemize}
	In particular, the category of complete filtered complexes has a model category structure. Because $\rmM^{\cfil}$ is enriched over $\strMod_k^{\cfil}$, it makes sense to talk about $\calP$-algebra objects in $\rmM^{\cfil}$. We will denote by $\strAlg_{\calP}(\rmM^{\cfil})$ the model category of $\calP$-algebras and $\Alg_{\calP}(\calC^{\cfil})$ the associated $\infty$-category. The same theorem works for bicolored operads.
\end{theorem}
\begin{example}\,
	\begin{itemize}
		\item All the classical operads $\bP_n$, $\Com$, $\Lie$, $\bP_n^{nu}$, $\Com^{nu}$ can be viewed in complete operads. They are defined, for instance, in \cite{LV}.
		\item The operad $\bP_n^{(nu)}$ is naturally filtered if we set the Poisson bracket to be of filtration degree $-1$, this operad will be denoted $\bP^{\gr(, nu)}_n$ and filtered algebras over that operad are $\bP^{(nu)}_n$ algebras such that the bracket decreases the filtration degree by $1$. The classical Hopf structure preserves filtration degrees, so it is a complete Hopf operad.
		\item The internal linear dual $\co\bP^{\gr}_n$ has comultiplication in filtration degree $0$ and cobracket in filtration degree $1$. It has a natural unital Hopf complete cooperad structure.
		\item The operad $\Lie^{\gr}$ of complete filtered $\Lie$ algebras. Those are algebra with a Lie bracket of filtration degree $-1$. This operad is isomorphic to $\Lie\{k(1)\}$.
	\end{itemize}
\end{example}

We will also need a filtered curved cooperad notion:
\begin{defn}
	A \defterm{curved complete cooperad} $\calD$ is a complete cooperad with a symmetric sequence morphism:
	\[\theta : \calD \rightarrow \mathbf{1}[2]\]
	satisfying
	\[(\theta\otimes \id_{\calD} - \id_{\calD}\otimes \theta)\circ \Delta = 0.\]
\end{defn}
There is an analog of deformation theory for curved cooperads given in \cite{curvkd}, with another kind of cobar construction, we give the construction here for completeness:
\begin{defn}{\cite[3.3.6.]{curvkd}}
	If $\calD$ is a curved coaugmented complete cooperad we let
	\[\Omega \calD \defeq \operatorname{Free}(\bar\calD[-1])\]
	with derivation $\d - \d_{\text{cobar}} + \d_{\theta}$ where $\d_{\theta}$ is the unique derivation extending 
	\[\bar\calD[-1] \overset{\theta}{\rightarrow} \mathbf{1}[1] \rightarrow \operatorname{Free}(\bar\calD[-1])[1].\]
	(and $\operatorname{Free}$ is the free operad). (It squares to $0$ according to \cite[3.3.7.]{curvkd})
\end{defn}
\begin{example}\,
\begin{itemize}
	\item The curved cooperad $\co\Lie^{\theta}$ is given by adjoining to each tree of $\co\Lie$ an extra composition by a curvature of degree $-2$. The deconcatenation preserves the positions of the curvature. The curvature of the cooperad is given by $\theta : \co\Lie^{\theta} \rightarrow \mathbf{1}[2]$ sending the two different compositions on the cobracket to $-1 * \id$. It is exactly the Koszul dual cooperad we get if we apply the unital Koszul duality developed in \cite{curvkd} for $\Com$. Hence, there is a quasi-isomorphism $\Omega\co\Lie^{\theta}\{1\} \eqrightarrow \Com$. 
	\item The curved cooperad $\co\bP_n^{\theta}$ built the same way as $\co\Lie^{\theta}$, there is a quasi-isomorphism $\Omega\co\bP_n^{\theta}\{k[n]\} \eqrightarrow \bP_n$, the operad of unital Poisson algebras.
	\item In the filtered setting, we can as well define $\co\bP_n^{\theta,\gr}$, although with $\theta$ in filtration degree $1$.
\end{itemize}
\end{example}
\subsection{Convolution algebras}\,
\par Let $C$ and $D$ be two complete filtered complexes. Then there is an internal Hom $\underline\Hom(C,D)$ which is naturally filtered. Let $\calP$ be a complete operad and $\calD$ a coaugmented complete cooperad. The coaugmentation $k \rightarrow \calD$ splits $\calD \simeq k\oplus \bar\calD$ as a symmetric sequence. We let
\[\Conv(\calD,\calP)\defeq \prod_{n\geq 0} \underline\Hom_{\Sigma_n}\left(\bar\calD(n),\calP(n)\right)\]
be the \defterm{Convolution Lie algebra} of $\calD$ and $\calP$. It is endowed with the usual Lie bracket that is compatible with filtrations. We will also need the non-reduced convolution algebra:
\[\Conv^0(\calD,\calP)\defeq \prod_{n\geq 0} \underline\Hom_{\Sigma_n}\left(\calD(n),\calP(n)\right).\]
We will need the following proposition that is stated in \cite[2.3]{curvlie}:
\begin{prop}
	There is a cobar construction $\Omega \calD$ for coaugmented cooperads $\calD$ and any morphism $\Omega \calD \rightarrow \calP$ to any operad $\calP$ gives rise to a Maurer--Cartan element in $\Conv\left(\calD,\calP\right)$.
\end{prop}
\par We have resolutions of operads:
\begin{prop}
	There are quasi-isomorphisms:
	\begin{itemize}
		\item $\Omega \co\Lie\{1\} \eqrightarrow \Com^{nu}$,
		\item $\Omega \co\Com^{nu}\{1\} \eqrightarrow \Lie$,
		\item $\Omega \co\bP_n^{\gr,\nu}\{k(-1)[n]\} \eqrightarrow \bP_n^{\gr,nu}$. Here the shift by $k(-1)$ is important because it will shift the comultiplication from filtration degree $0$ to filtration degree $-1$ which is what we want. And
		\item $\Omega\co\bP_n^{\gr,\theta}\{k(-1)[n]\}\eqrightarrow \bP_n^{\gr}$.
	\end{itemize}
\end{prop}
\begin{proof}
	The first two are trivially filtered, so it reduces to the known result for classical operads (see e.g. \cite[Chapters 7 and 8]{LV}). For the third one, the usual map $\Omega \co\bP_n\{n\} \rightarrow \bP_n$ is compatible with the filtration we have fixed. To check that it is a weak equivalence of complete filtered complete we have to check it is on the associated graded complex. Remark that once we forget the graded structure, we get the quasi-isomorphism $\Omega\co\bP_n\{n\}\rightarrow\bP_n$. But this forgetful functor is conservative, so the graded morphism is a quasi-isomorphism as well. For the last one, note that the curvature $\theta$ has filtration degree $0$ thanks to the shift.
\end{proof}
\section{Forms, formal geometry, and forms over formal prestacks}\label{sec:firststep}
\par In this section, we firstly recall the theory of derived prestacks as initiated in \cite{HAGII}. We then set the recent general definitions of shifted symplectic structures given historically in \cite{PTVV} and generalized in \cite{sympgrpd} for more general twists than shifts by $n$. All of this will be done in \cref{ssec:prest}. This setting is well adapted to talk about shifted symplectic structure of derived prestacks that are not necessarily Artin stacks. This is the case for formal derived prestacks firstly defined in \cite{GRII}. It was conjectured for a long time that the Lurie--Pridham Theorem (\cite[IV]{SAG},\cite{pridef}) relating deformation of the point to $\Lie$-algebras could be extended into a theorem relating (formal) deformation of derived schemes and Lie algebroids over it. It was proven in the first place in \cite{KozLie} that it was possible to do it for eventually coconnective finitely presented $k$-algebras $A$ (see \cite{surveyform} for a good survey on the subject), and then generalized for any coherent $k$-algebras $A$ in \cite{deffol}. Because it is not straightforward to explicitly construct Lie algebroids, and that derived foliations of \cite{fol} appear more naturally in this paper, we will need to relate Lie algebroids to derived foliations. This is the main purpose of \cite{follie}. We will recall the main theorems of those in \cref{ssec:def1}, \cref{ssec:def2}, and \cref{ssec:def3}. We then establish in \cref{ssec:formform} and \cref{ssec:formlag} the main theorem of this section, namely that if $\calF$ is a foliation over $A$, and $\calF \rightarrow A(0)$ has a Lagrangian structure of weight $-1$, then
\[\Spec(A)\rightarrow[\Spec(A)/\calF]\]
has a Lagrangian structure, where $[\Spec(A)/\calF]$ is the formal derived prestack classified by $\calF$.
\subsection{Prestacks, forms, and symplectic structures}\label{ssec:prest}
\par We relate the beginning of \cite{sympgrpd} here.
\par A \defterm{derived prestack} is an accessible functor $\CAlg_k^\conn \rightarrow \Sp$. We denote the category of prestacks by $\dPSt_k$. For $S \in \dPSt_k$ we denote by $\dPSt_S$ the category $\dPSt_{k/S}$. In the sequel, we will only name such object ``prestacks" by throwing away the adjective ``derived". For $A \in \CAlg_k^{\conn}$ we write $\Spec(A)$ for the functor corepresented by $A$.\\
Recall that if $X \in \dPSt_k$ we let
\[\QCoh(X) \defeq \lim_{\Spec(A) \rightarrow X} \Mod_A.\]
be the right Kan extension of the module category functor. We call this category the category of quasi-coherent sheaves over $X$. Those are symmetric monoidal categories, and if $f : X \rightarrow Y$ is a morphism of prestacks, we have an induced adjunction on quasi-coherent sheaves:
\[f^* : \QCoh(Y) \rightleftarrows \QCoh(X) : f_*\]
where $f^*$ is symmetric monoidal. We let $\Perf(X)$ be the subcategory of $\QCoh(X)$ of dualizable objects. Those are called \defterm{perfect complexes} over $X$. The following base prestacks will be of interest:
\begin{itemize}
	\item $S = B^{\pre}\bG_m$, the classifying prestack of $\bG_m$. It classifies graded complexes: $\QCoh(X \times B^{\pre}\bG_m) \simeq \QCoh(X)^{\gr}$.
	\item The quotient prestack $(\bA^1/\bG_m)^{\pre}$ classifying filtered complexes: $\QCoh(X\times (\bA^1/\bG_m)^{\pre}) \simeq \QCoh(X)^{\fil}$ (\cite{geomfil}).
\end{itemize}

If $\calA \in \CAlg(\QCoh(X))$ we have a relative spectrum construction:
\[\Spec(\calA)(f : \Spec(B)\rightarrow X) \defeq \Map_{\CAlg_B}(f^* \calA, B).\]
$\Spec(\calA)$ exists over $X$.
Recall from \cite[B.10.2.]{CHS} that it makes sense to ask when a morphism $f : X \rightarrow S$ admits a relative cotangent complex $\cot_{X/S}$. If it exists, it is unique, and we refer to it as the \defterm{cotangent complex} of $f$. When the cotangent complex exists as perfect complex, we denote by $\tan_X$ its dual, and we call it the \defterm{tangent complex} of $X$.
Let $\calL \in \QCoh(S)$ be an invertible object. Let $X \rightarrow S$ be a prestack over $S$. If $S = \Spec(A)$, we set
\[\DR_S(X)\defeq \lim_{\Spec(B) \rightarrow X \rightarrow S} \DR(A/B) \in \CAlg(\Mod_A^{\cfil}).\]
If $S$ is not affine anymore, we set 
\[\DR_S(X) \defeq \lim_{f : S'=\Spec(B) \rightarrow S} f_*\DR_{S'}(X\times_{S}S') \in \CAlg(\QCoh(S)^{\cfil}).\]
Recall from \cite{sympgrpd} and \cite{PTVV} that we can define the space of closed $\calL$-twisted $p$-form on $X$ to be
\[\calA^{p,\cl}(X/S,\calL) \defeq \Map_{\QCoh(S)^{\cfil}}\left(\mathbf{1}(p)\otimes[-p],\DR_S(X)\otimes\calL\right)\]
\begin{example}\,
	\begin{itemize}
		\item If $S = *$ and $\calL = k[n]$ we recover the classical notion of closed $n$-shifted $p$-forms.
		\item If $S = \bA^1/\bG_m^{\pre}$ and $\calL = k[n](-1)$ we will get the notion of closed $n$-shifted $p$-form of filtration degree $1$.
	\end{itemize}
\end{example}
The functor $X \mapsto \calA^{p,\cl}(X/S,\calL)$ defines a prestack over $S$.
\begin{defn}
	An \defterm{$\calL$-twisted pre-symplectic structure} on $X$ is the data of an $\calL$-twisted $2$-form on $X$. If $f : L \rightarrow X$ is a morphism of prestacks over $S$. An \defterm{$\calL$-twisted isotropic structure} on $f$ is the data of a commuting square
	\[\begin{tikzcd}
		L \ar[r,"f"] \ar[d]& X \ar[d,"\omega_X"]\\
		*\ar[r,"0"]	& \calA^{2,\cl}(\vari/S,\calL).
	\end{tikzcd}\]
	We denote by $\Isot(f/S,\calL)$ the space of $\calL$-twisted isotropic structures.
\end{defn}
\begin{remark}
	In \cite{sympgrpd} and \cite{PTVV}, they use mixed graded structures in order to define these structures. But it is fine since complete filtered objects and mixed graded objects form the same $\infty$-categories.
\end{remark}
\par In order to talk about non-degeneracy of pre-symplectic and isotropic structures, one has to compare those data with actual $p$-forms:
\begin{thm}\label{thm:grDRglobalform}\cite[2.6.]{sympgrpd}
	If $X\rightarrow S$ admits a cotangent complex then,
	\[\gr\left(\DR_S(X)\right) \simeq \Gamma(X,\Sym(\cot_{X/S}(1)[-1])).\]
\end{thm}
\begin{defn}
	Assume $X$ admits a perfect cotangent complex (relative over $S$). An \defterm{$\calL$-twisted symplectic structure} on $X$ is the data of a pre-symplectic structure $\omega$ on $X$ such that the induced morphism
	\[\Theta_{\omega} : \tan_{X/S} \rightarrow \cot_{X/S} \otimes \pi^*\calL\]
	is an equivalence ($\pi : X \rightarrow S$ is the canonical projection). We denote by $\Symp(X/S,\calL)$ the space of $\calL$-twisted symplectic structures.\\
	If $f : L \rightarrow X$ is a morphism admitting a perfect relative cotangent complex, an \defterm{$\calL$-twisted Lagrangian structure} on $f$ is the data of an isotropic structure $h$ on $f$, such that the pre-symplectic structure is symplectic and such that the induced morphism
	\[\Theta_h : \tan_f \rightarrow \cot_{L/S} \otimes \pi^* \calL[-1]\]
	is an equivalence. We denote by $\Lag(f/S,\calL)$ the space of $\calL$-twisted Lagrangian structures on $f$. When $\calF = k[n]$ we will often refer to such structures as ``$n$-shifted structures," and when $\calL = k(-1)[n]$, we will often refer to such structures as ``$n$-shifted structures of weight $-1$."
\end{defn}
When $S$ is the point, we omit all the occurence of ``$/S$" in the notations above.
\par Here, we have defined everything for prestacks, but note that we can define the same way internal twisted symplectic and Lagrangian structures on an object $A$ in $\CAlg(\calC)$. In that case we will keep the same notations but replace $\Spec(A)$ by $A$ and $\Spec(f)$ by $f$ for $f$ a morphism of commutative rings. For instance, $\Symp(A,\mathbf{1}[n])$ will denote the space of $n$-shifted symplectic structure on $A$.
\subsection{Finiteness properties on algebras}\label{ssec:def1}
\par Let $A \in\CAlg^\conn_k$. Let us briefly summarize finiteness properties on $A$ and modules over $A$ in order to properly do formal geometry under $\Spec(A)$.
\begin{defn}\cite[7.2.4.16]{HA}
	$A$ is called \defterm{coherent} if:
	\begin{itemize}
		\item Every finitely generated ideal of $\pi_0(A)$ is finitely presented,
		\item $\pi_i(A)$ are finitely presented $\pi_0(A)$-modules for $i>0$.
	\end{itemize}
\end{defn}
This assumption will be relevant because the theory of pro-coherent sheaves is well-behaved over such rings \cite{deffol}.
\begin{defn}\cite{HA}
	The category $\Mod_A$ of $A$-modules is endowed with a t-structure $(\Mod_A^{\geq 0},\Mod_A^{\leq 0})$. An $A$-module $M$ is called 
	\begin{itemize}
		\item \defterm{perfect} if it is a dualizable object of $\Mod_A$,
		\item \defterm{almost perfect} if there exists for each $n$ an $n$-connective map $P_n \rightarrow M$ where $P_n$ is perfect, and
		\item \defterm{coherent} if it is almost perfect and eventually coconnective.
	\end{itemize}
	We let $\Perf_A$, and $ \APerf_A \supset \Coh_A$ be the full subcategories of $\Mod_A$ spanned by, respectively, perfect, almost perfect, and coherent objects.
\end{defn}

\begin{defn}
	A morphism of (filtered, graded) connective algebras $B\rightarrow A$ is said to be \defterm{almost of finite presentation} if for any $n$, $\tau_{\leq n}A$ is a compact object in the category $\tau_{\leq n}\CAlg_B$. The algebra $A$ is said to be \defterm{almost of finite presentation} if $k\rightarrow A$ is almost of finite presentation. Remark that, because $k$ is a field, almost of finite presentation algebras are coherent.
\end{defn}
\begin{lm}\cite[A.19.]{deffol}
	If $B \rightarrow A$ is almost of finite presentation, then $\cot_{A/B}$ is almost perfect.
\end{lm}
\par We end this section by giving the following definition:
\begin{defn}
	A map $B\rightarrow A$ of connective commutative $k$-algebras is a \defterm{nilpotent extension} if $\pi_0(B)\rightarrow \pi_0(A)$ is surjective with nilpotent kernel.
\end{defn}
\subsection{Lie algebroids and derived foliations}\label{ssec:def2}
\par Assume now that $A$ is coherent. We define
\[\QCoh_A^{\vee}\]
the category of \defterm{pro-coherent modules} over $A$ to be the category of functors
\[\Mod_A^{\conn} \rightarrow \Sp\]
preserving
\begin{itemize}
	\item pullbacks along $\pi_0$-surjections,
	\item the terminal object,
	\item limits of almost eventually constant towers (towers that become eventually constant when applying $\tau_{\leq i}$), and
	\item filtered colimits of $n$-coconnective objects for any $n$.
\end{itemize}
By virtue of \cite[A.42.]{deffol}, we have an equivalence $\QCoh_A^{\vee} \simeq \Ind(\Coh_A^\op)$. By virtue of \cite[A.28.]{deffol} it is stable and presentable, and by of \cite[A.45.]{deffol} it is presentably symmetric monoidal. Because $A$ is a $k$-algebra, it is also $k$-linear.
\begin{remark}
	We can see that the definition here recovers exactly the object $\mathrm{Pro}(\QCoh^-(\Spec(A)))_{\mathrm{laft}}$ that is used in \cite[Chapter 1, 3.4.1.]{GRII}.
\end{remark}
\begin{warn}\label{warn:pro}
	Beware that in \cite{deffol}, the authors use \textit{homological} convention whereas in \cite{GRII}, the authors use the \textit{cohomological} convention. Also, in \cite{deffol}, the authors write $\Ind(\calC^\op)$ for the category of pro-coherent objects of $\calC$, but there is a little ambiguity as $\mathrm{Pro}(\calC)$ should be $\Ind(\calC^\op)^\op$. They have chosen to throw away the extra $\op$ for simplicity. It doesn't really matter since functoriality of $\QCoh^\vee_{(\vari)}$ will remain the same if we put it or not.
\end{warn}
For $X$ a prestack, we set
\[\QCoh^\vee(X) \defeq \lim_{\Spec(A)\rightarrow X} \QCoh_A^\vee.\]
There is a fully faithful duality functor:
\[\APerf_A^\op\overset{(\vari)^{\vee}}{\longrightarrow} \QCoh_A^{\vee}.\]
The essential image will be denoted by $\APerf_A^{\vee}$, its objects are called \defterm{dually almost perfect complexes} over $A$.
\par Recall from \cite[3.8.]{deffol} that there is a sifted colimit-preserving monad $\Lie_A$ defined on $\QCoh_A^{\vee}$ corresponding to \defterm{partition Lie algebroids} over $A$. We let $\LieAlgd_A \defeq \Alg_{\Lie_A}(\QCoh^{\vee}_A)$ be the category of \defterm{Lie algebroids} over $A$ (we omit the substantive ``partition" because we work in characteristic zero). There is a ``Koszul duality" functor:
\[\calD : \Art_A^{\op} \rightarrow \LieAlgd_A\]
coming from a deformation context in the sense of \cite[IV]{SAG}. The image spans $\LieAlgd_A$ by sifted colimits.\\
\par Recall as well that Fu defines in \cite[4.25.]{follie} a sifted colimit-preserving monad called $S$ on $\QCoh_A^{\vee,\fil,\leq 0}$ intuitively corresponding to Lie algebroids with a bracket of filtration degree $-1$. Then Koszul duality provides a diagram
\[\begin{tikzcd}
	\Alg_S(\QCoh_A^{\vee,\fil,\leq 0}) \ar[r,shift left, "\CE^{\fil}"] \ar[d,shift left,"(\vari)^u"]& \ar[l, shift left, "\calD^{\fil}"] \CAlg_{k,\gr(\vari)^0 \simeq A}^{\fil,\geq 0,\op} \ar[d, shift left, "\aug"]\\
	\LieAlgd_A \ar[u, shift left, "\const"] \ar[r, shift left, "\CE"] & \ar[l, shift left, "\calD"] \ar[u, shift left, "\adic"] \CAlg_{k/A}^{\op},
\end{tikzcd}\]
where $\CAlg_{k,\gr(\vari)^0 \simeq A}^{\fil,\geq 0}$ is the category consisting of an object $B$ of $\CAlg(\QCoh_k^{\fil,\geq 0})$, and an equivalence $\gr(B)^0 \simeq A$. Note that it is not exactly what Fu writes. In \cite{follie}, it is more natural to work with pro-coherent modules because the author works over an arbitrary base ring. But in our setting it is irrelevant since the base ring, $k$, is a field, and we have that canonically $\Mod_k^{\vee} \simeq \Mod_k$. The commuting left adjoints are $\aug \circ \CE^{\fil} \simeq \CE\circ(\vari)^u$. The filtered Chevalley--Eilenberg functor is then defined to be
\[\cCE \defeq \CE^{\fil}\circ \const.\]
By virtue of \cite[4.27.]{follie}, $\cCE$ actually ends in $\CAlg^{\cfil,\geq 0,\op}$ so the notation is consistent.
\begin{prop}
	$\cCE$ preserves sifted colimits.
\end{prop}
\begin{proof}
	$\CE^{\fil}$ sends colimits to limits since it is a left adjoint. By definition, the monads $\Lie_B$ and $S$ are sifted colimit-preserving monads. Since $\const$ preserves sifted colimits on the level of modules, it preserves them also on the level of algebras. Remark that the condition $\gr(\vari)^0 \simeq A$ is harmless with respect to limits.
\end{proof}
\begin{defn}
	The category $\Fol$ of \defterm{derived foliations} is defined to be the pullback of the following diagram:
	\[\begin{tikzcd}
		\Fol \ar[r,"\iota"] \ar[d,"e"] & \CAlg^{\cfil} \ar[d,"\gr"]\\
		\Mod \ar[r,"\bigwedge^\bullet"] & \CAlg^{\gr}.
		\arrow["\lrcorner"{anchor=center, pos=0.125}, draw=none, from=1-1, to=2-2]
	\end{tikzcd}\]
	Where $\Mod$ is the category where objects are $(B,M)$ where $B \in \CAlg$ and $M \in \Mod_B$ and $\bigwedge^\bullet$ sends $(B,M)$ to $\Sym(M(1)[-1])$. We set $\Fol_A \defeq \Fol\times_{\CAlg}\{A\}$. For $\calF \in \Fol$, we will use the following notation: $e(\calF) = (\gr^0\calF,\cot_{\calF})$. The complexThe complex  $\cot_{\calF}$ is called the \defterm{cotangent complex} of $\calF$. Because in the sequel, we will also consider $\cot_{\calF}$ computed internally to complete filtered objects, and because the new notation is in competition with this one, we will instead denote by $\cot_\calF^\int$ the cotangent complex computed \textit{internally} to complete filtered objects. We denote by $\tan_\calF$ the linear dual of $\cot_\calF$ and $\tan_\calF^\int$ the internal linear dual of $\cot_\calF^\int$.
\end{defn}
\begin{example}
	According to \cref{lm:grdr}, if $B \rightarrow A$ is a morphism of commutative $k$-algebras, then $\DR(A/B)$ is a derived foliation over $A$. The cotangent complex $\cot_{\DR(A/B)}$ of the foliation is actually $\cot_{A/B}$ so we usually write the second one instead of the first one. We still refer to the cotangent complex computed internally to complete filtered objects as $\cot_{\DR(A/B)}^\int$.
\end{example}
This is the combination of theorem \cite[4.26.]{follie} and the remark \cite[5.12.]{follie}:
\begin{thm}\label{thm:lie=fol}
	$\cCE$ induces an equivalence:
	\[\LieAlgd_A^{\dap} \simeq \Fol^{\ap,\op}_A.\]
	where $\LieAlgd_A^{\dap}$ are the Lie algebroids with underlying dually almost perfect $A$-module and $\Fol^{\ap}_A$ are the derived foliations with underlying almost perfect $A$-module.
\end{thm}
\begin{lm}
	If $B \in \Art_A$ then $\cCE(\calD(B)) \simeq \DR(A/B)$.
\end{lm}
\begin{proof}
	It follows from the diagram defining $\cCE$ that we naturally get an equivalence $\aug(\cCE(\calD(B))) \simeq (B \rightarrow A)$ so we get an arrow $\adic(B\rightarrow A) \rightarrow \cCE(\calD(B))$ which is an equivalence on the associated graded (by definition the underlying object of $\calD(B)$ is $\tan_{A/B}$). Since everything is complete, it is an equivalence ($\adic(B\rightarrow A) \simeq \DR(A/B)$).
\end{proof}
\par If we combine the fact that $\cCE$ preserves sifted colimits, $\LieAlgd_A$ is spanned by sifted colimits of objects in the image of $\calD$, and this lemma, we get:
\begin{cor}\label{cor:drspansfol}
	Every $\calF \in \Fol_A^{\ap}$ is a cosifted limit of foliations of the form $\DR(A/B)$ for $B \rightarrow A$ Artin. The underlying cotangent complex is also the cosifted limit of the $\cot_{A/B}$.
\end{cor}
\subsection{Formal geometry and formal integration of derived foliations}\label{ssec:def3}
In this section, we recall the definition of formal prestacks and recall the main theorem of \cite{deffol}.
\begin{defn}
	A map $f : X \rightarrow Y$ of prestacks is said to be:
	\begin{itemize}
		\item \defterm{convergent} if for any $A \in \CAlg_k^{\conn}$, the following diagram:
			\[\begin{tikzcd}
				X(A) \ar[r] \ar[d] & \lim X(\tau_{\leq n} A)\ar[d]\\
				Y(A) \ar[r] & \lim Y(\tau_{\leq n} A)
			\end{tikzcd}\]
			is cartesian. We say that $X$ is \defterm{convergent} if $X\rightarrow *$ is convergent. We denote by $\dPSt_k^{\conv}$ the full subcategory of convergent prestacks. It is a localization of $\dPSt_k$ with left adjoint $X \mapsto \lim X(\tau_{\leq n}(\vari))$.
		\item \defterm{locally almost finitely presented} or \defterm{of finite presentation} if it is convergent, and if for any $n\geq 0$ and any filtered diagram $I \rightarrow \CAlg_k^{\conn}$ consisting only of $n$-truncated rings and of colimit $A$, the following square
			\[\begin{tikzcd}
				\colim X(A_{\alpha}) \ar[d] \ar[r] & X(A)\ar[d]\\
				\colim Y(A_{\alpha}) \ar[r] & Y(A)
			\end{tikzcd}\]
			is cartesian. The prestack $X$ is \defterm{locally almost of finite presentation} if $X \rightarrow *$ is locally almost finitely presented.
	\end{itemize}
\end{defn}
\begin{remark}\label{rmk:lafpconv}
	If $X\rightarrow *$ satisfies only the second condition of the definition, then its \defterm{convergent completion} $X^{\conv} \defeq \lim X(\tau_{\leq n}(\vari))$ is locally almost of finite presentation. Because a filtered colimit of $n$-truncated rings is always $n$-truncated.
\end{remark}
\begin{warn}
	In \cite{SAG}, Lurie refers to ``convergent" prestacks as ``nilcomplete" prestacks. We have chosen the denomination of \cite{GRI}.
\end{warn}
\par If $A \in \CAlg_k$, we set $\Ind\Coh_A \defeq \Ind(\Coh(A))$ and
\[\Ind\Coh(X) \defeq \lim_{\Spec(A)\rightarrow X} \Ind\Coh_A.\]
The category $\Ind\Coh(X)$ is called the category of \defterm{ind-coherent sheaves} over $X$.
Ind-coherent sheaves are relevant because they are really well suited to do formal geometry. For locally almost of finite presentation prestacks, it has a natural symmetric monoidal structure and an action of quasi-coherent sheaves on it. It also has a nice duality property that allows us to define a functor:
\[\Upsilon_X : \QCoh(X) \rightarrow \Ind\Coh(X).\]
This functor is natural with respects to pullbacks, that is, for $f : X \rightarrow Y$, there is a commutative square:
\[\begin{tikzcd}
	\QCoh(Y) \ar[r,"\Upsilon_Y"] \ar[d,"f^*"] & \Ind\Coh(Y)\ar[d,"f^!"]\\
	\QCoh(X) \ar[r,"\Upsilon_X"] & \Ind\Coh(X).
\end{tikzcd}\]
We will need the following (that is also defined in \cite[B.18.]{deffol}):
\begin{defn}
	A prestack $X$ is said to \defterm{admit a pro-coherent tangent complex} if for any $x : \Spec(B) \rightarrow X$, the functor
	\[M \in \Mod_B^{\conn} \mapsto \Map_{\Spec(B)/}(\Spec(B[M]),X)\]
	defines a pro-coherent sheaf on $B$ and if the family $(\cot_{X,x})_{x\in X(B)}$ defines an object of $\QCoh^\vee(X)$. In that case, we denote by $\tan_X$ the complex, by remembering that it lives in $\QCoh^\vee(X)$.
\end{defn}
\begin{warn}
	In \cite{GRII}, the authors talk about pro-coherent cotangent complexes, the definition is the same here and the only difference comes from the difference between the category of pro-coherent sheaves we consider (see \cref{warn:pro}) (because the functoriality of the tangent/cotangent complex will be different).
\end{warn}
\par Serre duality (\cite[Chapter 5, 4.2.10.]{GRI}) also gives us:
\begin{lm}\label{lm:pro=ind}
	If $B$ is an almost of finite presentation $k$-algebra, then there is an equivalence.
	\[\QCoh^{\vee}_B \simeq \Ind\Coh_B\]
\end{lm}
\par Under this equivalence, the functor interwines the duality functor $\Upsilon_X$ with the functor
\[\iota_X : \begin{array}{ccc}\QCoh(X) &\longrightarrow& \QCoh^\vee(X)\\ M & \longmapsto & M\otimes_{\calO_X} (\vari) \end{array}\]
If $f : B \rightarrow A$, there is a functor
\[\QCoh_B^\vee \rightarrow \QCoh_A^\vee\]
that is given by precomposition by
\[f_* : \Mod_A^+ \rightarrow \Mod_B^+\]
where $\Mod_A^+$ denotes the full subcategory of $\Mod_A$ of eventually coconnective modules (see \cite[A.37.]{deffol}). This functor is a pullback functor and $\iota_{(\vari)}$ is natural w.r.t. this pullback. We will call $f^!$ this functor because under Serre duality we also recover Gaitsgory--Rozenblyum's shriek pullback functor. This functor is also symmetric monoidal (\cite[A.45.]{deffol}).
\begin{warn}
	Here, we use $\Ind\Coh(X)$ for convenience because every statement in the sequel comes from \cite{GRI, GRII}, that use $\Ind\Coh(X)$, but we believe that they should also work if we replace $\Ind\Coh(X)$ by $\QCoh^\vee(X)$ without using Serre duality and by only assuming that $X$ is locally coherent.
\end{warn}
\begin{lm}\label{lm:perfindcoh}
	If $X$ is locally almost of finite presentation, then the functor mentioned above induces a fully faithful embedding
	\[\Perf(X) \rightarrow \Ind\Coh(X)\]
	with essential image dualizable objects. If $f : X \rightarrow Y$ is a map between locally almost of finite presentation prestacks, then
	\[f^! : \Ind\Coh(Y) \rightarrow \Ind\Coh(X)\]
	preserves dualizable objects.
\end{lm}
\begin{proof}
	The first part is \cite[Chapter 6, 3.3.7.]{GRI}. For the second part, the argument is similar. Locally, for $S \defeq \Spec(B) \rightarrow X$ with $B$ almost of finite presentation and eventually coconnective, the image of $\calE$ in $\Ind\Coh(S)$ is dualizable, but because $\QCoh(S)$ fully faithfully embeds in $\Ind\Coh(S)$, it is perfect in $\QCoh(S)$, it is then perfect in $\QCoh(X)$.
\end{proof}

\begin{defn}
	We set $X_{\dR}$ to be the prestack defined by $A \mapsto X(A^{\red})$. It is called the \defterm{de Rham} prestack of $X$. We also set $X_{\red}\defeq \colim_{\Spec(A) \rightarrow X} \Spec(A^{\red})$ to be the image along the left Kan extension of the reduction functor on $X$.\\
	If $f : X \rightarrow Y$ is a map of prestack, $Y^{\wedge}_X \defeq Y \times_{Y_{\dR}} X_{\dR}$ is called the \defterm{formal completion} of $X$ along $f$.
\end{defn}
\begin{defn}
	A prestack $X$ is said to \defterm{have a deformation theory} if $X \rightarrow *$ is convergent and if for any cartesian square
	\[\begin{tikzcd}
		B' \ar[r]\ar[d] & B\ar[d]\\
		A'\ar[r] & A
	\end{tikzcd}\]
	with $B \rightarrow A$ a nilpotent extension, the induced square
	\[\begin{tikzcd}
		X(B') \ar[r]\ar[d] & X(B)\ar[d]\\
		X(A')\ar[r] & X(A)
	\end{tikzcd}\]
	is cartesian.
\end{defn}
\par We now have every definition at hands to define formal maps.
\begin{defn}
	A map $X \rightarrow Y$ of prestacks is a \defterm{formal thickening}  if:
	\begin{itemize}
		\item $Y$ has deformation theory,
		\item $X \rightarrow Y$ is locally almost finitely presented,
		\item and $X \rightarrow Y$ induces an equivalence $X_{\red} \eqrightarrow Y_{\red}$.
	\end{itemize}
	We let $\Thick(X) \hookrightarrow \dPSt_{X/}$ be the category of formal thickenings under $X$. We say that $Y$ is a \defterm{formal thickening} of $X$ when there exists $f : X \rightarrow Y$ that is a formal thickening.
\end{defn}
\begin{remark}
	In \cite{deffol}, the authors refer to formal thickenings as ``formal moduli stacks". If $X$ and $Y$ are moreover supposed to be locally almost of finite type, we get the notion of ``formal moduli" under $X$ of \cite{GRII}.
\end{remark}
\begin{defn}
	If $X$ has a perfect cotangent complex, we define the space $\LagThick(X,n+1)$ of \defterm{$(n+1)$-shifted Lagrangian thickenings} of $X$ to be the full subgroupoid of the fiber product
	\[\Isot\Thick(X,n) \defeq \Thick(X) \times_{\dPSt_{k, X/}} \Isot\Map(X, n+1)\]
	consisting of isotropic structures that are Lagrangian, where $\Isot\Map(X, n+1)$ is the category $\Fun(\Delta^2,\dPSt_k)\times_{\Fun(\Delta^1 \vee \Delta^1,\dPSt_k)} \{X\rightarrow * \overset0\rightarrow \calA^{2,\cl}(\vari,n+1)\}$.
\end{defn}
\par Formal thickenings under an affine scheme $\Spec(A)$ can be more explicitly described; they are completely determined by Artin extensions of $A$:
\begin{defn}
	An \defterm{Artin extension} of $A$ is a map $B\rightarrow A$ of connective commutative $k$-algebras that satisfies
	\begin{itemize}
		\item $B \rightarrow A$ is a nilpotent extension,
		\item $B \rightarrow A$ is almost of finite presentation,
		\item the fiber is eventually coconnective.
	\end{itemize}
	We let $\Art_A$ be the full subcategory of $\CAlg_{k/A}$ of Artin extensions of $A$.
\end{defn}
\begin{defn}
	A \defterm{formal moduli problem} under $\Spec(A)$ is an accessible functor
	\[F : \Art_A \rightarrow \Sp\]
	satisfying the following conditions:
	\begin{itemize}
		\item $X(A) \simeq *$,
		\item For any cartesian square in $\Art_A$
		\[\begin{tikzcd}
			B \ar[r] \ar[d] & B_0\ar[d]\\
			B_1\ar[r] & B_{01},
		\end{tikzcd}\]
		such that $B_0 \rightarrow B_{01}$ is a nilpotent extension, the induced square
		\[\begin{tikzcd}
			X(B) \ar[r]\ar[d] & X(B_0) \ar[d]\\
			X(B_1) \ar[r] & X(B_{01})
		\end{tikzcd}\]
		is cartesian.
	\end{itemize}
	We let $\Moduli_A$ be the category of formal moduli problems under $\Spec(A)$. The Yoneda embedding defines a fully faithful functor $\Art_A^{\op} \rightarrow \Moduli_A$.
\end{defn}
The main theorem \cite[4.19.]{deffol} is:
\begin{thm}\label{thm:liealgdfmp}
	If $A$ is coherent, there is an equivalence
	\[\Moduli_A \simeq \LieAlgd_A.\]
\end{thm}

Thanks to \cite[section 5.1.]{deffol}, it is possible to compare formal moduli problems with formal thickenings the following way:
\begin{thm}\label{thm:fmpthick}
	There is an equivalence $\Thick(\Spec(A)) \eqrightarrow \Moduli_A$ sending a formal thickenings $\Spec(A) \rightarrow Y$ to the functor defined on each Artin extension $B \rightarrow A$ by the space of lifts
	\[\begin{tikzcd}
		\Spec(A) \ar[r] \ar[d] & Y\\
		\Spec(B). \ar[ur,dashed]
	\end{tikzcd}\]
	The fact that affines are sent to affines implies that the diagram
	\[\begin{tikzcd}
		\Moduli_A \ar[r] & \dPSt_k\\
		\Art_A^{\op} \ar[u]\ar[r] & \CAlg_k^{\conn,\op}\ar[u]
	\end{tikzcd}\]
	commutes.
\end{thm}
\par If we combine \cref{thm:lie=fol,thm:liealgdfmp,thm:fmpthick} we obtain:
\begin{prop}
	If $A$ is coherent, there is a fully faithful embedding
	\[\Fol_A^{\ap,\op} \rightarrow \Thick(\Spec(A)).\]
	If $\calF$ is such a derived foliation and $X\defeq\Spec(A)$, we denote by $[X/\calF]$ the associated formal prestack. We call it the associated \defterm{formal leaves prestack}. 
\end{prop}
\begin{prop}\label{prop:folspnice}
	If $X = \Spec(A)$ is almost of finite presentation and $\calF\in\Fol_A^{\ap}$ with $\pi : X\rightarrow [X/\calF]$ the projection. Then $[X/\calF]$ satisfies the following:
	\begin{itemize}
		\item It is locally almost of finite presentation.
		\item It admits a pro-coherent tangent complex.
		\item We always have a fiber sequence of pro-coherent complexes:
			\[\tan_\calF \rightarrow \tan_X \rightarrow \pi^!\tan_{[X/\calF]}.\]
	\end{itemize}
\end{prop}
\begin{proof}We prove each statement separatly:
	\begin{itemize}
		\item For the first assertion, it is known thanks to \cite[section 5.1.]{deffol} that $[X/\calF]$ is a convergent completion of a left Kan extension along the inclusion $[\Art_A,\Sp] \hookrightarrow \dPSt_{X_{\dR}}$. That is, before taking the completion, it is a colimit consisting of objects of the form $\Spec(B)$ for $B \in \Art_A$ (they are almost of finite presentation according to \cref{lm:artinafp}), or of the form $\Spec(A)_{\dR}$ which is locally almost of finite presentation. Hence it satisfies the second condition of the definition of locally almost of finite presentation prestacks and $[X/\calF]$ is locally almost of finite presentation by \cref{rmk:lafpconv}.
		\item For the second assertion, $[X/\calF]$ has a pro-coherent tangent complex because it has a deformation theory (\cite[Chapter 1, 7.2.5.]{GRII}) and because it is locally almost of finite presentation (\cite[Chapter 1, 3.5.2.]{GRII}).
		\item For the last assertion, recall that according to \cite[section B.3.]{deffol}, one has to consider the almost finitely presented approximation of derivation functors in order to control the functoriality of the tangent complexes. But here we don't have to do it since the cotangent complexe exists as pro-coherent objects (\cite[A.26.]{deffol}). Then, thanks to \cite[5.26.]{deffol} and \cite[B.28.]{deffol}, there is the fiber sequence
			\[\tan_\calF \rightarrow \tan_X \rightarrow \pi^!\tan_{[X/\calF]}.\]
	\end{itemize}
\end{proof}
\begin{remark}\label{rmk:intfibseq}
	We have a similar sequence if we compute internal cotangent complexes. If $p$ denote the projection map $\calF \rightarrow A(0)$, observe that $\cot_p \simeq \cot_{\calF}(1)$. This is because of the fiber sequence:
	\[p^* \cot_\calF^\int \rightarrow \cot_A \rightarrow \cot_p.\]
	If we compute the associated graded, we get a fiber sequence
	\[\cot_A\oplus \cot_{\calF}(1)[-1] \rightarrow \cot_A \rightarrow \gr\cot_p\]
	giving an equivalence $\gr\cot_p \simeq \cot_{\calF}(1)$, which lifts in the complete filtered world since it is fully concentrated in one weight. Finally, we get naturally a fiber sequence:
	\[|p^*\cot_\calF^\int|\rightarrow \cot_A \rightarrow \cot_{\calF}\]
	implying the following identity:
	\[|p^*\cot_\calF^\int| \simeq \pi^!\cot_{[X/\calF]}.\]
\end{remark}
\begin{prop}\label{prop:folperf}
	If $X = \Spec(A)$ is of finite presentation, and $\calF$ is a derived foliation over $A$ with a perfect cotangent complex, then $[X/\calF]$ has a perfect cotangent complex and we have a fiber sequence
	\[\pi^*\cot_{[X/\calF]} \rightarrow \cot_X \rightarrow \cot_{\calF}.\]
	Moreover, $\cot_\calF^\int$ is perfect as a complete filtered $\calF$-module.
\end{prop}
\begin{proof}
	Thanks to \cref{prop:folspnice}, $[X/\calF]$ is locally almost of finite presentation, so we can apply the results of \cite[Chapter 5]{GRII}. If $\pi: \Spec(A) \rightarrow [\Spec(A)/\calF]$ denotes the projection, we have a fiber sequence of ind-coherent complexes:
	\[\tan_\calF \rightarrow \tan_X \rightarrow \pi^!\tan_{[X/\calF]}\]
	Then $\pi^!\tan_{[X/\calF]}$ is perfect. But by \cite[Chapter 5, 2.2.6.]{GRII}, $\pi^!$ factors through an equivalence:
	\[\Tot\Ind\Coh(R_{\bullet}) \simeq \Ind\Coh([X/\calF])\]
	where $R_{\bullet}$ is a groupoid over $X$. Because the exceptional pullbacks preserves dualizable objects by \cref{lm:perfindcoh}, and by \cite[4.6.1.11.]{HA}, $\pi^!\tan_{[X/\calF]}$ is dualizable in $\Tot\Ind\Coh(R_{\bullet})$ so it is dualizable in $\Ind\Coh([X/\calF])$. This implies that $[X/\calF]$ has a perfect cotangent complex, and that we have a fiber sequence:
	\[\pi^*\cot_{[X/\calF]} \rightarrow \cot_X \rightarrow \cot_{\calF}.\]
	For the second part, it is because $\gr$ reflects dualizable objects, which is also a direct consequence of \cref{lm:grhom}. It is then enough to check that $\cot_\calF^\int$ is dualizable on the associated graded, which is clear by the assumptions.
\end{proof}
\subsection{Forms on formal leaves prestacks}\label{ssec:formform}
\begin{lm}\label{lm:DRconv}
	Let $X \rightarrow Y$ be a morphism of prestacks inducing an equivalence on the associated convergent prestacks. Then $\DR(X) \rightarrow \DR(Y)$ is an equivalence.
\end{lm}
\begin{proof}
	This is because $\DR(\vari)$ is convergent by \cite[3.1.]{cohohall3CY}. Being an equivalence on the associated convergent prestacks means that both stacks have the same eventually coconnective affine sites, so the de Rham algebras are the same.
\end{proof}
\begin{prop}\label{prop:DRFMP}
	The de Rham functor $\DR: \Moduli_A^\op \rightarrow \CAlg_{k/\DR(A)}^{\cfil}$ obtained by composing the usual one on prestacks with the inclusion $\Moduli_A \rightarrow \dPSt_{\Spec{A}/}$ is right Kan extended from the inclusion $\Art_A^{\op} \rightarrow \Moduli_A$ and the de Rham functor $\Art_A \rightarrow \CAlg_{k/A} \rightarrow \CAlg_{k/\DR(A)}^{\cfil}$.
\end{prop}
\begin{proof}
	The functor $\Moduli_A \rightarrow \dPSt$ comes from a series of colimit-preserving functors:
	\[[\Art_A,\Sp] \rightarrow [\mathrm{Ext}_A,\Sp] \rightarrow \dPSt^{\conv}_{\Spec(A)//\Spec(A)_{\dR}}\]
	where $\mathrm{Ext}_A$ is the category of Artin-like extensions of $A$ but without the almost finite presentability assumption. 
	The functor $\DR: \dPSt \rightarrow \CAlg_k^{\cfil,\op}$ is the unique colimit-preserving functor extending $\DR: \CAlg_k^{\op} \rightarrow \CAlg_k^{\cfil,\op}$. It extends into a colimit-preserving functor $\dPSt_{\Spec(A)/} \rightarrow \CAlg_{k/\DR(A)}^{\cfil,\op}$. By \cref{lm:DRconv} and \cite[5.5.4.20.]{HTT}, it factors as a colimit-preserving functor $\dPSt^{\conv}_{\Spec(A)/} \rightarrow \CAlg_{k/\DR(A)}^{\cfil,\op}$. Because $[\Art_A,\Sp]$ is the free cocompletion of $\Art_A^{\op}$, and $[\Art_A,\Sp] \rightarrow \dPSt_{\Spec(A)/}$ sends representable functors to affine derived schemes, $\DR: [\Art_A,\Sp] \rightarrow \CAlg_{k/\DR(A)}^{\cfil,\op}$ is obtained as the left Kan extension of $\DR: \Art_A^{\op} \rightarrow \CAlg_{k/A}^{\op} \rightarrow \CAlg_{k/\DR(A)}^{\cfil,\op}$. Now, since $\Moduli_A$ is a full subcategory of $[\Art_A,\Sp]$, we have the result by transitivity of Kan extensions. 
\end{proof}
\begin{defn}
	We let $\DR^{\falg}: \Fol_A \rightarrow \CAlg_k^{\cfil}$ be the ``formally algebraic" de Rham functor given by
	\[\DR^{\falg}\defeq |\DR^{\int} \circ \iota|^{\int}\]
	where $\iota: \Fol_A \rightarrow \CAlg_k^{\cfil}$ is the inclusion. The realization is the internal realization. To be more precise: $\DR^{\int}(A)$, where $A$ is a complete filtered commutative algebra, has two filtrations. One comes from the filtration of $A$ (internal filtration), and the other one comes from the weight of differential forms (external filtration). The functor $|\vari|^{\int}$ applies $|\vari|$ to each component of the external filtration. Because we will deal with two filtrations in the sequel, we will refer to $M\langle 1 \rangle$ as $M$ put in external degree/weight $1$ and $M(1)$ as $M$ put in internal degree/weight $1$.
\end{defn}
We will prove the following theorem:
\begin{thm}\label{thm:algform}
	Let $\calF \in \Fol_A^{\ap}$ with $A$ almost of finite presentation, then there is a natural equivalence:
	\[\DR([X/\calF]) \simeq \DR^{\falg}(\calF).\]
	It is built in such a way that there is a commuting diagram in $\Mod_k^{\gr}$:
	\[\begin{tikzcd}
		\left|\Sym_{\calF}(\cot_\calF^\int\langle 1 \rangle[-1])\right|^{\int} \ar[d]&\Gamma\left([X/\calF],\Sym_{\calO_{[X/\calF]}}(\cot_{[X/\calF]}\langle 1 \rangle[-1])\right)\ar[d]\ar[l,"\sim"']\\
		\left|\Sym_{A(0)}(\cot_\calF^\int\otimes_{\calF} A(0)\langle 1 \rangle[-1])\right|^{\int}&\Gamma\left(X,\Sym_{\calO_{X}}(\pi^*\cot_{[X/\calF]}\langle 1 \rangle[-1])\right)\ar[l,"\sim"']
	\end{tikzcd}\]
	The first horizontal map comes from the first equivalence, \cref{lm:grdr}, and \cref{thm:grDRglobalform}.
\end{thm}
Several lemmas will be needed.
\begin{lm}\label{lm:artinafp}
	If $f: B \rightarrow A$ is an Artin extension with $A$ almost of finite presentation, then $B$ is almost of finite presentation.
\end{lm}
\begin{proof}
	Let $M$ be the fiber of $f$. We have a long exact sequence:
	\[\dots \rightarrow \pi_1(A) \rightarrow \pi_0(M) \rightarrow \pi_0(B) \rightarrow \pi_0(A) \rightarrow 0.\]
	The module $\pi_0(M)$ has a finitely presented nilpotent image. Since $\pi_0(A)$ is finitely presented, $\pi_0(B)$ is also finitely presented. By \cite[A.19.(4)]{deffol}, $M$ is almost perfect as a $B$-module, so $\pi_i(M)$ is a finite type as $\pi_0(B)$-module; hence, $\pi_i(B)$ are of finite type as $\pi_0(B)$-modules. Hence, $B$ is almost of finite presentation.
\end{proof}
\begin{lm}\label{lm:adiccotap}
	Assume $B \rightarrow A$ is Artin with $B$ and $A$ almost of finite presentation. Then $\adic(B \rightarrow A)$ is complete and almost of finite presentation, and $\cot_{\adic(A \rightarrow B)}$ is almost perfect.
\end{lm}
\begin{proof}
	By \cite[2.18.]{deffol} we know that such a map satisfies the condition $(3)$ of \cite[2.15.]{deffol}. So that $\adic(A \rightarrow B) \rightarrow B$ is complete almost finite augmented in the sense of \cite[2.15.(2)]{deffol}. Because the extension is Artin, by \cref{lm:artinafp}, $(k\rightarrow k) \rightarrow (B \rightarrow A)$ is almost of finite presentation. Now, $\adic$ preserves almost of finite presentation objects (see the proof of \cite[2.18.]{deffol}). So $\adic(B \rightarrow A)$ is almost of finite presentation. This implies by \cite[A.19.(2).]{deffol} that $\cot_{\adic(B \rightarrow A)}$ is almost perfect as a $\adic(B \rightarrow A)$-module.
\end{proof}

\begin{lm}\label{lm:algformaff}
	Let $B\rightarrow A$ be an Artin extension, with $A$ almost of finite presentation. We have:
	\[\DR(B) \simeq |\DR^{\int}(\DR(A/B))|^{\int}.\]
\end{lm}
\begin{proof}
	First, observe that the associated graded on the right-hand side is
	\[|\Sym_{\DR(A/B)}(\cot_{\DR(A/B)}^\int(1)[-1])|^{\int}.\]
	By evaluating at $0$ we get $|\DR(A/B)| \simeq F^0(\adic(B \rightarrow A)) \simeq B$ (\cref{lm:adicfil}), so it induces a map
	\[\DR(B) \rightarrow |\DR^{\int}(\DR(A/B))|^{\int}.\]
	In order to check that it is an equivalence, we check it on the associated graded. On each component, the map is the following:
	\begin{equation}\label{grintdr}\left(\cot_B[-1]^{\otimes_B n}\right)_{h\Sigma_n}\rightarrow \left|\left(\cot_{\DR(A/B)}^\int[-1]^{\otimes_{\DR(A/B)} n}\right)_{h\Sigma_n}\right|\end{equation}
	Now remark that since $B \rightarrow A$ is Artin with $A$ almost of finite presentation, by \cref{lm:adiccotap}, $\cot_{\adic(B\rightarrow A)}$ is almost perfect. Note also that since $F^0(\adic(B\rightarrow A))\simeq B$ and $\cot_{B/A}$ is $1$-connective ($\pi_0(B) \rightarrow \pi_0(A)$ is surjective so the cotangent is $1$-connective by \cite[7.4.3.2.]{HA}), we have that $\adic(B\rightarrow A)$ is connective for the neutral t-structure on filtered objects. By \cite[A.22.]{deffol}, $\cot_{\adic(B\rightarrow A)}^{\otimes_A n}$ is complete for any $n$ since for any $n$ those are also almost perfect. Since they are all concentrated in non-negative filtration degree, we can perform the computation of the right-hand side of \eqref{grintdr} in filtered objects instead of complete one. Since the underlying filtered object is symmetric monoidal and a left adjoint, it commutes with everything in the equation. We conclude because $\left(\DR(A/B)\right)^u \simeq \left|\DR(A/B)\right| \simeq B$.
\end{proof}
\begin{lm}\label{lm:formaldesc}
	Assume $\calF \in \Fol_A^{\ap}$ is a cosifted limit of $\calF_i \in \Fol_A^{\ap}$ and that $\cot_{\calF}$ is also the cosifted limit of $\cot_{\calF_i}$. Then the natural functor:
	\[\DR^{\falg}(\calF) \rightarrow \lim \DR^{\falg}(\calF_i)\]
	is an equivalence.
\end{lm}
\begin{proof}
	Before computing the natural morphism, remark that $\gr \calF_i \simeq \Sym_A(\cot_{\calF_i}(1)[-1])$ and that $\lim \cot_{\calF_i} \simeq \cot_{\calF}$. Remark that the natural morphism is induced by the realization of the natural morphism
	\begin{equation}\label{eqdesc}\DR^{\int}(\calF) \rightarrow \lim \DR^{\int}(\calF_i)\end{equation}
	because $|\vari|$ commutes with limits. So it is enough to check this is an equivalence. By looking at the associated graded, we have
	\[\Sym_{ \calF}(\cot_\calF^\int\langle 1\rangle[-1]) \rightarrow \lim\Sym_{\calF_i}(\cot_{\calF_i}^\int\langle 1\rangle[-1]).\]
	We can in the same way consider the associated graded on each component and check it is an equivalence here. Since taking the associated graded commute with limits and the $\Sym$ functor, we have that the left-hand side is
	\[\Sym_{\gr\calF}(\cot_A \otimes_A \gr\calF\langle 1 \rangle[-1] \oplus \cot_{\calF} \otimes_A \gr\calF\langle 1 \rangle(1)[-1]),\]
	and the right-hand side is
	\[\lim \Sym_{\gr\calF_i}(\cot_A \otimes_A \gr\calF_i\langle 1 \rangle[-1] \oplus \cot_{\calF_i} \otimes_A \gr\calF_i\langle 1 \rangle(1)[-1]).\]
	We can simplify, and we get
	\[\Sym_{A}(\cot_A \langle 1 \rangle[-1] \oplus \cot_{\calF}\langle 1 \rangle(1)[-1] \oplus \cot_{\calF}(1)[-2]) \rightarrow \lim \Sym_{A}(\cot_A \langle 1 \rangle[-1] \oplus \cot_{\calF_i}\langle 1 \rangle(1)[-1] \oplus \cot_{\calF_i}(1)[-2]).\]
	This is an equivalence by a graded and sifted version of \cite[3.34.]{follie}, and because we are in zero characteristic. Because $\gr$ is conservative, the map \eqref{eqdesc} is an equivalence.
\end{proof}
\begin{lm}\label{lm:tensorreal}
	If $A$ is almost of finite presentation and $B \rightarrow A$ is an Artin extension
	\[\left|\cot_{\DR(B/A)}^\int\otimes_{\DR(B/A)} A(0)\right|\simeq \left|\cot_{\DR(B/A)}^\int\right|\otimes_{\left|\DR(B/A)\right|} A.\]
\end{lm}
\begin{proof}
	This is an analog of \cref{lm:algformaff}.
\end{proof}
\begin{proof}[Proof of \cref{thm:algform}]
	If $\calF \in \Fol_A^{\ap}$, by \cref{cor:drspansfol}, it is a cosifted limit of $\DR(A/B)$. So if we apply \cref{lm:formaldesc}, \cref{prop:DRFMP}, and \cref{lm:algformaff} it proves the theorem. For the commuting diagram, the reasoning is similar. First, we have that the following diagram:
	\[\begin{tikzcd}
		\left|\cot_{\DR(A/B)}^\int\right|^{\otimes n}_{h\Sigma_n} \ar[d] & {\cot_{|\DR(A/B)|}}^{\otimes n}_{h\Sigma_n}\ar[l]\ar[d]\\
		\left|\cot_{\DR(A/B)}^{\int,\otimes n}\otimes_{\DR(A/B)}A(0)\right|_{h\Sigma_n} & \left(\cot_{|\DR(A/B)|}^{\otimes n}\otimes_{|\DR(A/B)|} A\right)_{h\Sigma_n}\ar[l]
	\end{tikzcd}\]
	commutes. By the same argument as in the proof of \cref{lm:algformaff}, the upper horizontal map is an equivalence, and the upper right corner is equivalent to ${\cot_B}^{\otimes_n}_{h\Sigma_n}$. By \cref{lm:tensorreal}, the bottom horizontal map is an equivalence, and the bottom right corner is equivalent to $\left(\cot_B^{\otimes n}\otimes_B A\right)_{h\Sigma_n}$. Hence we get a commutative diagram:
	\[\begin{tikzcd}
		\left|\Sym_{\DR(A/B)}(\cot_{\DR(A/B)}^\int\langle 1 \rangle[-1])\right|^{\int} \ar[d] & \Sym_B(\cot_B\langle 1 \rangle[-1])\ar[l,"\sim"']\ar[d]\\
		\left|\Sym_{A(0)}(\cot_{\DR(A/B)}^\int\otimes_{\DR(A/B)} A(0)\langle 1 \rangle[-1])\right|^{\int} & \Sym_A(\cot_B \otimes_B A\langle 1 \rangle[-1]).\ar[l,"\sim"']
	\end{tikzcd}\]
	Now we can apply a similar descent argument. Thanks to the two fiber sequences
	\[\cot_{[X/\calF]}\otimes_{[X/\calF]} A \rightarrow \cot_A \rightarrow \cot_{\calF}\]
	and
	\[|\cot_\calF^\int\otimes_{\calF} A(0)| \rightarrow \cot_A \rightarrow \cot_{\calF},\]
	we know that they are both limits of, respectively, $\cot_B\otimes_B A$ and $|\cot_{\DR(A/B)}^\int\otimes_{\DR(A/B)}A(0)|$ (recall that $\cot_{\calF}$ is a limit of $\cot_{A/B}$). We conclude because $\Sym$ commutes with cosifted limits of almost perfect objects (\cite[3.34.]{follie}).
\end{proof}
\begin{remark}
	\cref{thm:algform} remains true in a slightly more general setting where we allow to integrate any Lie algebroids instead of dually almost perfect one. One has to find the correct internal de Rham analog for Lie algebroids. It turns out that the $\gr^0$ functor which is lax monoidal induces a map from complete filtered Lie algebroids (where the Lie structure is of filtration degree $0$) to Lie algebroids, and it has a right adjoint $(\vari)^{\dR}$. This right adjoint commutes with sifted colimits so we can adapt the argument of \cref{lm:formaldesc} in order to state that $\DR(X_{\calL}) \simeq \CE(\calL^{\dR})$ for a Lie algebroid $\calL$ with associated formal moduli problem $X_{\calL}$. The algebra $\CE(\calL^{\dR})$ is classically called the \defterm{Weil algebra} of $\calL$ (see e.g. \cite{weilalg}). A computation of forms can be found in \cite{sympliealgd}.
\end{remark}
\subsection{Lagrangian structures on formal leaves prestacks}\label{ssec:formlag}
\par We are now able to compare algebraic Lagrangian structures on a derived foliation with Lagrangian structures on its formal leaves prestack.
\begin{thm}\label{thm:lagalglagthick}
	Let $X \defeq \Spec(A)$ almost of finite presentation and $\calF$ a derived foliation with almost perfect cotangent complex on $A$. There is a canonical morphism $p: \calF \rightarrow A(0)$ inducing the morphism $\pi: X \rightarrow [X/\calF]$. There is an equivalence
	\[\Isot(\pi,k[n]) \simeq \Isot(p,k[n](-1)).\]
	According to \cref{prop:folperf}, if $A$ is of finite presentation and $\calF$ has a perfect cotangent complex, then $[X/\calF]$ has a perfect cotangent complex, and $\calF$ has a perfect internal cotangent complex. It then makes sense to talk about symplectic and Lagrangian structures. The former equivalence fits into a commutative square:
	\[\begin{tikzcd}
		\Isot(\pi,k[n]) \ar[r,"\sim"] & \Isot(p,k[n](-1))\\
		\Lag(\pi,k[n]) \ar[r,"\sim"] \ar[u,hook] & \Lag(p,k[n](-1)). \ar[u,hook]
	\end{tikzcd}\]
\end{thm}
\begin{proof}
	For convenience, we will write $(\vari)\langle n \rangle$ for $(\vari)(n)$ coming from the external filtration degree (the one that appears when we compute the de Rham complex), and we do not change the notation for the internal one (the one coming from the foliation structure). We also write $\calC \defeq \Mod_k^{\cfil}$, and we denote by $\DR^{\int}(\vari)$ the de Rham functor internal to $\calC$. Before proving the equivalence, let's remark the following fact: Let $C$ be the cofiber
	\[\DR^{\int}(\calF)\otimes k(-1) \rightarrow \DR^{\int}(\calF) \rightarrow C.\]
	Since $|\vari |^{\int}$ is exact, it induces a fiber sequence
	\[|\DR^{\int}(\calF)\otimes k(-1)|^{\int} \rightarrow |\DR^{\int}(\calF)|^{\int} \rightarrow |C|^{\int}\]
	implying that $|C|^{\int} \simeq \gr^{0,\int}(\DR^{\int}(\calF)) \simeq \DR(\gr^0(\calF)) \simeq \DR(A)$ by \cref{prop:drcomtrivgr}.\\
	We now have:
	\begin{align*}	\Isot(p,k[n](-1)) & \simeq \Map_{\calC^{\cfil}}(k\langle 2 \rangle[-n-2](0),\DR^{\int}(\calF)\otimes k(-1))\times_{\Map_{\calC^{\cfil}}(k\langle 2 \rangle[-n-2](0),\DR^\int(A(0))\otimes k(-1))} *\\
		&\simeq \Map_{\calC^{\cfil}}(k\langle 2 \rangle[-n-2](0),\DR^{\int}(\calF)\otimes k(-1))\\
		&\simeq \Map_{\calC^{\cfil}}(k\langle 2 \rangle[-n-2](0),\DR^{\int}(\calF))\times_{\Map_{\calC^{\cfil}}(k\langle 2 \rangle[-n-2](0),C)} *\\
		&\simeq \Map_{\Mod_k^{\cfil}}(k\langle 2 \rangle[-n-2],|\DR^{\int}(\calF)|^{\int})\times_{\Map_{\Mod_k^{\cfil}}(k\langle 2 \rangle[-n-2],|C|^{\int})}*\\
		&\simeq \Map_{\Mod_k^{\cfil}}(k\langle 2 \rangle[-n-2],\DR([X/\calF]))\times_{\Map_{\Mod_k^{\cfil}}(k\langle 2 \rangle[-n-2],\DR(A))}*\\
		&\simeq \Isot(\pi,k[n]).
	\end{align*}
	The second equivalence is because $\Map_{\calC^{\cfil}}(k\langle 2 \rangle[-n-2](0),\DR^\int(A(0))\otimes k(-1))$ is contractible, essentially because there is no non-trivial map from $k(0)$ to $k(-1)$. Now let's check the non-degeneracy conditions.
	Assume that $X$ and $\calF$ have a perfect cotangent complex. Let $(\omega,h)$ a $k(-1)[n]$-twisted isotropic structure on $p$ with associated isotropic structure $(\omega',h')$ on $\pi$. We have the following diagram:
	\[\begin{tikzcd}
		p^*\tan_\calF^\int(1)[-n] \ar[r] \ar[d,"p^*\Theta_\omega"]&\tan_p(1)[-n+1]\ar[r]\ar[d,"\Theta_h"]&\tan_{A(0)}(1)[-n+1]\ar[d,"(-1)^{-n+1}\Theta_h^\vee(1){[-n+1]}"]\\
		p^*\cot_\calF^\int \ar[r] & \cot_{A(0)} \ar[r] & \cot_p.
	\end{tikzcd}\]
	If we consider the realization we obtain
	\[\begin{tikzcd}
		\pi^*\tan_{[X/\calF]}[-n] \ar[r]\ar[d,"\pi^*\Theta_{\omega'}"] & \tan_{\calF}[-n+1] \ar[r]\ar[d,"\Theta_{h'}"] & \tan_A[-n+1] \ar[d,"(-1)^{-n+1}\Theta_{h'}^\vee{[-n+1]}"]\\
		\pi^*\cot_{[X/\calF]} \ar[r] & \cot_A \ar[r] & \cot_{\calF}.
	\end{tikzcd}\]
	$|\tan_p(1)|\simeq \tan_\calF$ and $|\cot_p|\simeq \cot_\calF$ by the remark \cref{rmk:intfibseq}. And the pullback on $A$ of the induced form $\omega'$ is the induced form of the pullback on $A(0)$ thanks to the commuting diagram of \cref{thm:algform}.
	The vertical arrows in the upper diagram are equivalences if and only if the vertical arrows in the bottom diagram are equivalences. We conclude because:
	\begin{itemize}
		\item $\pi^*$ reflects equivalences between perfect objects, by a similar reasoning as in \cref{prop:folperf}.
		\item $p^*$ reflects this equivalence, by looking at the associated graded.
	\end{itemize}
	This means that the structure on $p$ is Lagrangian if and only if the structure on $\pi$ is Lagrangian.
\end{proof}

\begin{remark}
	The $C$ in the proof is a complete filtered object with the underlying object being $\DR(A)$, we can think that it provides a non-trivial extra filtration on $\DR(A)$, but if we remark that it is given by the ``internal associated graded" (every filtered object is internally filtered when we tensor by $\dots \rightarrow k(-1) \rightarrow k(0) \rightarrow k(1) \rightarrow \dots$) which is the associated graded seen as a filtered object with $0$ transition maps, so the extra filtration is trivial in that sense.
\end{remark}	

\section{Complete filtered Poisson and coisotropic structures}\label{sec:sndstep}
\par In this section, we give a recollection of the results of \cite{MSI}, but adapted to the complete filtered setting. In \cref{ssec:defop1}, \cref{ssec:defop2}, \cref{ssec:defop3}, and \cref{ssec:defop4} we mostly cover section $3$ of \cite{MSI}. In \cref{ssec:poisstruct} and \cref{ssec:coisstruct} we will cover section $4$ of \cite{MSI}. Recall from \cref{sec:modelcat} that everything happening in this section will take place in a model category $\rmM$. In particular, algebras will be taken strictly in the model category of algebra objects of $\rmM$.
\subsection{Filtered braces}\label{ssec:defop1}
Let $\calD$ be a complete filtered Hopf cooperad. Classically, if we forget the filtration, there is an operad $\Br_{\calD}$ called the \defterm{Brace operad} constructed in \cite{braces}. One can check that because $\calD$ is filtered, the operad $\Br_{\calD}$ is naturally filtered where a decorated tree is filtered by the sum of the filtrations of its decorations. Because we have asked all the structures to be compatible with the filtration and that the filtration is complete (on each arity, the filtration comes from the decoration, which is a finite product of element of a complete filtered complex), it is a complete operad. Let's recall conceptually how the brace construction work. First, given a (complete) operad $\calT$ with a morphism $\Lie \rightarrow \calT$, and a $\calT$-algebra $A$, it is possible to make sense of the notion of Maurer--Cartan element in $A$ (it is naturally a $\Lie$-algebra). One can then ask if, given a Maurer--Cartan element $\pi$, the twisted complex $A^\pi$ is still a $\calT$-algebra. This is not the case in general; however, it is possible to build a new operad $\mathrm{Tw}\calT$ such that $A^\pi$ is a $\mathrm{Tw}\calT$-algebra. This procedure is called the operadic twisting procedure and is described in \cite{optwist}. Now, given a (complete) operad $\calP$, it is possible to extend the $\preLie$-structure of
\[\Conv(\calD,\calP)\]
into a $\preLie_{\calD}$-algebra structure ($\preLie_\calD$ is an enhancement of the operad $\preLie$ that uses the Hopf structure of $\calD$). Given a Maurer--Cartan element $\pi$ of it,
\[\Conv(\calD,\calP)^\pi\]
have a natural $\mathrm{Tw}\preLie_{\calD}$-algebra structure. This is by definition the brace operad $\Br_\calD$ and this example is the main example of $\Br_{\calD}$-algebras that we will use. The same phenomenon happens for $\Conv^0(\calD,\calP)$ and a lot of convolution-like constructions. The most important fact about the brace construction is that its property really look like a kind of BV tensor product $\mathbb{E}_1 \otimes_{\mathrm{BV}} \calP$ of operads. In \cite{Padd}, Safronov builds a functor
\[\Alg_{\Br_\calD}(\calC) \rightarrow \Alg(\Alg_{\Omega \calD}(\calC))\]
and proves that in the case $\Omega \calD$ is a (resolution) of the operad $\bP_n$ of Poisson algebras it is an equivalence. This, combined with the following theorem of \cite{braces}, establishes Poisson additivity, an analog of Lurie--Dunn additivity for Poisson operads.
\begin{theorem}
	There is a quasi-isomorphism:
	\[\Omega\left(\co\bP_{n+1}^{\gr,\theta}\{1\}\right) \eqrightarrow \Br_{\co\bP_n^{\gr,\theta}}.\]
\end{theorem}
\begin{proof}
	For the same reason as the proof for Koszul duality between $\co\bP_n^{\gr,nu}$ and $\bP_n^{\gr,nu}$, we can check easily that the morphism built by Calaque and Willwacher in \cite{braces} is compatible with the filtrations, and that it is a weak equivalence of complete filtered complexes.
\end{proof}
\begin{defn}
	Let $\calD$ be a complete counital Hopf cooperad and $A$ a $\Omega\left(\calD\{\calL\}\right)$-algebra. Then the algebra structure is given by a Maurer--Cartan element $f$ of $\Conv(\calD\{\calL\};A)$ (\cite[2.29]{curvlie}). We set
	\[\Z(A)\defeq \Conv^0(\calD^{\cu}\{\calL\};A)^{f}\otimes\calL^{-1}\]
	twisted by the Maurer--Cartan $f$. It is called the \defterm{center} of $A$. It has a natural $\Br_{\calD}\{\calL\}$-algebra structure. Remark that it has a natural grading coming from the arities. The brace construction will be compatible with this external grading.
\end{defn}
\subsection{Relative braces}\label{ssec:defop2}
We can define as well the relative brace construction. Given two complete filtered complexes $A$ and $B$ and a complete cooperad $\calD$, we set
\[\Conv^0\left(\calD;A,B\right) \defeq \Hom\left(\calD(A),B\right).\]
Then it is possible to endow the triple
\[\left(\Conv^0\left(\calD^{\cu};A\right),\Conv^0\left(\calD;A,B\right),\Conv^0\left(\calD^{\cu};B\right)\right)\]
into a colored algebra over a ``relative $\preLie_{\calD}$-algebra" defined in \cite[3.2]{MSI}, the only modification we have to do is to allow $\preLie_{\calD}^{\rightarrow}$ to be filtered, with the filtration of an element being the sum of the filtrations of the decorations coming from $\calD$. This allows us to do relative deformation theory in the following sense:
\begin{prop}
	If $(C_1,C_2,C_3)$ is a complete filtered $\preLie^{\rightarrow}$-algebra then
	\[C_1 \oplus C_2[-1] \oplus C_3\]
	has a natural $\calL_{\infty}$-algebra structure (in complete filtered complexes).
\end{prop}
If it were also a $\preLie^{\rightarrow}_{\calD}$-algebra, given a Maurer--Cartan element $(f_1,f_2,f_3)$ in the obtained homotopy $\Lie$ algebra, we can relatively twist the complexes by it. The twisted complexes
\[\left(C_1^{f_1},C_2^{f_2},C_3^{f_3}\right)\]
will have a ``relative $\Br_{\calD}$-algebra" structure, by filtering $\Br^{\rightarrow}_{\calD}$ with the filtration coming from $\calD$. This allows us to do relative operadic twisting. Some new brace algebra structures arise that way.
\subsection{Swiss cheese construction}\label{ssec:defop3}
We need to adapt a little the Swiss-Cheese construction of \cite[3.3]{MSI}. We fix an invertible complete filtered complex $\calL$ instead of only a number. Let $\calD_1$ and $\calD_2$ be two coaugmented cooperads where the latter also has a Hopf counital structure. And let
\[F: \Omega\calD_1 \rightarrow \Br_{\calD_2}\{\calL\}\]
be a morphism of operads. Let $\{\calA,\calB\}$ be a set of colors. We define the $\{\calA,\calB\}$-colored operad $\SC(\calD_1,\calD_2)$ to be the free colored operad spanned by the colored symmetric sequence on $P(\calD_1,\calD_2)$ given by:
\begin{itemize}
	\item $P(\calD_1,\calD_2)(\calA^{\otimes n},\calA) \defeq \bar\calD_1(n)$,
	\item $P(\calD_1,\calD_2)(\calB^{\otimes l},\calB) \defeq \bar\calD_2\{\calL\}(l)$,
	\item $P(\calD_1,\calD_2)(\calA^{\otimes m}\otimes\calB^{\otimes l},\calB) \defeq \calD_1(m)\otimes \calD_2^{\cu}\{\calL\}(l)\otimes \calL[1]$.
\end{itemize}
with the differential given in \cite[3.12]{MSI} which are well defined since they are also compatible with filtration degrees.\\
Concretely, an $\SC(\calD_1,\calD_2)$-algebra is the data of a $\Omega\calD_1$-algebra $A$ and a $\Omega\calD_2$-algebra $B$ with an $\infty$-morphism of $\Omega\calD_1$-algebras $A \rightarrow \Z(B)$ ($\Z(B)$ is a $\Br_{\calD_2}\{\calL\}$-algebra hence a $\Omega\calD_1$-algebra).\\
One can explicitly construct the $\calL_{\infty}$-algebra classifying the deformation theory of $\SC(\calD_1,\calD_2)$-algebras. Let $\calL(\calD_1,\calD_2;A,B)$ be the $\calL_{\infty}$-algebra given by
\[\calL(\calD_1,\calD_2;A,B) \defeq \Conv(\calD_1;A)\oplus\Conv(\calD_2\{\calL\};B) \oplus \Hom(\calD_1(A)\otimes\calD_2^{\cu}\{\calL\}(B),B)\otimes\calL^{-1}[-1]\]
with brackets defined as in \cite[30]{MSI}. We have the following statement analogous to \cite[3.14.]{MSI}:
\begin{prop}
	\[\Map_{2\Op}(\SC(\calD_1,\calD_2),\End_{A,B}) \simeq \mathrm{MC}(\calL(\calD_1,\calD_2;A,B)).\]
\end{prop}

\subsection{Poisson structures}\label{ssec:poisstruct}
Let $A\in\rmM^{\cfil}$. We can endow $A$ with 
\begin{itemize}
	\item An external grading where $A$ is in weight $0$.
	\item A graded $\bP_{n+1}^{\gr}$-algebra structure with zero bracket.
\end{itemize}
\begin{defn}
	The complete filtered complex of \textbf{weak $n$-shifted polyvector fields} on $A$ is defined to be the center of $A$:
	\[\wPol^{\eps}(A,n) \defeq |\Z(A)|^{\cfil}\]
	It has a homotopy (externally, with a bracket of weight $-1$) graded $\bP^{\gr}_{n+2}$-structure (coming from the brace construction).
\end{defn}
We can also define the strict polyvectors:
\begin{defn}
	The complete filtered complex of \defterm{strict $n$-shifted polyvector fields} on $A$ is defined to be
	\[\Pol^{\eps}(A,n) \defeq \Hom_A\left(\Sym_A(|\Omega^1_A|^{\cfil}\langle -1 \rangle[n+1](-1)),A\right)\]
	where $|\Omega^1_A|^{\cfil}$ is the internal complete filtered complex of Kähler differentials. It is externally graded if we put $\Omega^1_A$ in weight $-1$, and it has a (externally) graded $\bP^{\gr}_{n+2}$-algebra structure coming from the Schouten--Nijenhuis bracket (it indeed decreases the filtration degree by $1$). We may later consider internal polyvectors
	\[\Pol^{\eps,\int}(A,n) \defeq \Hom_A\left(\Sym_A(\Omega^1_A\langle -1 \rangle[n+1](-1)),A\right).\]
	also with the Schouten--Nijenhuis bracket.
\end{defn}
\begin{prop}\label{prop:strictpoly}
	Let $A$ be a bifibrant complete filtered cdga. Then
	\[\wPol^{\eps}(A,n) \simeq \Pol^{\eps}(A,n)\]
	as (externally) graded homotopy $\bP_{n+2}^{\gr}$-algebras.
\end{prop}
\begin{proof}
	For readability, we denote $\calL \defeq k(-1)[n+1]$. We have:
	\begin{align*}
		\wPol(A,n) &\simeq \left|\Conv^0(\co\bP^{\gr,\theta,\cu}_{n+1}\{\calL\langle-1\rangle\};A)\otimes \calL^{-1}\langle 1\rangle\right|^{\cfil}\\
			   &\simeq \left|\Hom(\co\Com^{\cu}\left(\co\Lie^{\theta}\{\calL^{-1}[1]\langle 1\rangle\}(A\otimes \calL\langle-1\rangle)\right),A)\right|^{\cfil}\\
			   &\simeq \left|\Hom(\Sym(\co\Lie^{\theta}(A[1])\otimes \calL[-1]\langle-1\rangle),A)\right|^{\cfil}\\
			   &\simeq \left|\Hom_{A}(\Sym_A(\co\Lie^{\theta}(A[1]) \otimes A \otimes \calL[-1]\langle-1\rangle),A)\right|^{\cfil}\\
			   &\simeq \left|\Hom_{A}(\Sym_A(\left|\mathrm{Harr}_{\bullet}(A,A)\right|^l\otimes \calL\langle-1\rangle),A)\right|^{\cfil}\\
			   &\simeq \left|\Hom_{A}(\Sym_A(\Omega^1_A\otimes \calL\langle-1\rangle),A)\right|^{\cfil}
			   \simeq \Pol^{\eps}(A,n).
	\end{align*}
	Here, we use the formula $\calD\{\calL\}(A) \simeq \calD(A\otimes\calL)\otimes\calL^{-1}$ and $\mathrm{Harr}(A,A)\simeq \Omega^1_A$ (complete filtered analog of \cite[2.6]{MSI} that uses the cofibrancy condition). Fibrancy condition is useful because in that case $\Hom_A(\vari,A)$ preserves weak equivalences.
\end{proof}
When $A \in \calC^{\cfil}$, we can define mixed graded shifted Poisson structures on $A$:
\begin{defn}
	The space $\Pois^{\eps}(A,n)$ of \defterm{complete filtered $n$-shifted Poisson structures} on $A$ is the maximal subgroupoid of the fiber of
	\[\Alg_{\bP^{\gr}_{n+1}}(\calC^{\cfil}) \rightarrow \CAlg(\calC^{\cfil})\]
	at $A$.
\end{defn}
\begin{prop}\label{prop:bivec}
	There is an equivalence:
	\[\Pois^{\eps}(A,n)\simeq \Map_{\Alg_{\Lie\{k\langle 1\rangle(1)\}}(\Mod_k^{\gr,\cfil})}\left(k\langle 2\rangle[-1],\wPol^{\eps}(A,n)\otimes k(-1)[n+1] \right).\]
\end{prop}
\begin{proof}
	Because we have settled up the deformation theory of complete filtered operads, the proof is the complete filtered analog of \cite[4.5]{MSI}.
\end{proof}
Thanks to this theorem, we can now talk about non-degeneracy of Poisson structures:
\par Given a filtered $n$-Poisson structure on $A$ and $\tilde A$ a bifibrant resolution of $A$, we have an $\infty$-morphism of graded $\Lie$-algebras $k(2)[-1] \rightarrow \Pol^{\eps}(\tilde A,n)\otimes k(-1)[n+1]$ by \cref{prop:bivec} and \cref{prop:strictpoly}. Which in turn gives a Maurer--Cartan element in $\Pol^{\eps}(\tilde A,n)\otimes k(-1)[n+1]$ by \cite[1.19]{MSI} (when we forget the filtration). But $\pi$ is in filtration degree $0$ by definition, thus $[\pi,\vari]$ defines a map:
\[\theta_\pi : \Omega^1_{\tilde A}\rightarrow T_{\tilde A}[-n](1)\]
(where $T_{\tilde A}$ is the tangent sheaf of $\tilde A$). We then get a map
\[\theta_\pi: \cot_A \rightarrow \tan_A[-n](1)\]
of complete filtered $A$-modules.
\begin{defn}
	We say that $\pi$ is \defterm{non-degenerate} if $\theta_\pi$ is an equivalence. We let $\Pois^{\eps,nd}(A,n)$ be the maximal subgroupoid of $\Pois^{\eps}(A,n)$ spanned by non-degenerate Poisson structures.
\end{defn}
Note that non-degeneracy can be checked on the underlying associated graded since it is conservative and computations of cotangent and tangent commute with taking the associated graded.\\
One can relate non-degenerate shifted Poisson structures to shifted symplectic structures:
\begin{thm}\label{thm:poisndsymp}
	There is an equivalence:
	\[\Pois^{\eps,nd}(A,n) \simeq \Symp(A,k(-1)[n]).\]
\end{thm}
\begin{proof}
	This is the same proof as in \cite{CPTVV} or \cite{poisndpri}; this is a consequence of \cref{prop:cois2pois} and \cref{thm:coisndlag} that has a proof in \cref{sec:proofnd}.
\end{proof}

\subsection{Relative Poisson algebras}\label{ssec:defop4}
We now introduce the filtered version of the operad $\bP_{[n+1,n]}$. Let $A$ be a strict $\bP^{\gr}_n$-algebra. The Poisson bracket on $A$ gives strictly a bivector $\pi$ of degree $1-n$ and of filtration degree $-1$ on $A$ satisfying the strict Maurer--Cartan equation, we can define the strict Poisson center of $A$:
\[\rZ(A) \defeq \Hom_A(\Sym(\Omega^1_A\langle -1 \rangle[n](-1)),A)\]
along with the Schouten--Nijenhuis bracket and the differential twisted by $[\pi,\vari]$.
\begin{defn}
	A $\bP^{\gr}_{[n+1,n]}$-algebra consists of:
	\begin{itemize}
		\item A $\bP^{\gr}_{n+1}$-algebra $A$,
		\item a $\bP^{\gr}_n$-algebra $B$,
		\item and a morphism of $\bP^{\gr}_{n+1}$-algebras $A \rightarrow \rZ(B)$.
	\end{itemize}
	There is a forgetful functor $\Alg_{\bP^{\gr}_{[n+1,n]}}(\calC^{\cfil}) \rightarrow \left(\CAlg(\calC^{\cfil})\right)^{\Delta^1}$ sending $(A,B,A\rightarrow \rZ(B))$ to $A \rightarrow \rZ(B) \rightarrow B$.
\end{defn}
The Swiss-Chesse construction of earlier allows us to explicitly build a cofibrant resolution of this operad. We set
\[\tilde\bP^{\gr}_{[n+1,n]} \defeq \SC(\co\bP^{\gr,\theta}_{n+1}\{k(-1)[n+1]\},\co\bP^{\gr,\theta}_n).\]
It is a cofibrant complete operad because it is quasi-free (by definition). There is a natural morphism:
\[\tilde\bP^{\gr}_{[n+1,n]} \rightarrow \bP^{\gr}_{[n+1,n]}\]
coming from the map $\Z(B) \rightarrow \rZ(B)$ built exactly as in the proof of \cref{prop:strictpoly} (and the fact that strict Poisson structures are a special case of homotopy Poisson structures). One can prove exactly as in \cite[3.19.]{MSI} the following:
\begin{prop}
	The natural morphism $\tilde\bP^{\gr}_{[n+1,n]} \rightarrow \bP^{\gr}_{[n+1,n]}$ is a weak equivalence.
\end{prop}
\subsection{Coisotropic structures}\label{ssec:coisstruct}
Let $f: A \rightarrow B$ be a morphism of commutative algebras in $\calC^{\cfil}$. Following \cite[4.2.]{MSI} we can define the complete filtered complex of weak relative $k(-1)[n]$-twisted polyvector fields $\wPol^{\eps}(f,n)$ that is a homotopy $\bP^{\gr,nu}_{n+2}$-algebra. Again, it has a strict model denoted $\Pol^{\eps}(f,n)$ that is equivalent (compatible only a priori with the underlying $\calL_{\infty}$ structure) to the weak one when $f$ is bifibrant and the natural morphism $\Pol(A,n) \rightarrow \Pol_A(B,n-1)$ is surjective. Weak relative polyvectors are compatible with weak polyvectors so that we have a diagram:
\[\begin{tikzcd}
	&\wPol^{\eps}(f,n)^{\geq 2}[n+1]\ar[rd]\ar[ld]&\\
	\wPol^{\eps}(A,n)^{\geq 2}[n+1] && \wPol^{\eps}(B,n-1)^{\geq 2}[n].
\end{tikzcd}\]
of $\calL_{\infty}$-algebras.

\begin{defn}
	The space $\Cois^{\eps}(f,n)$ of \defterm{complete filtered $n$-shifted coisotropic structures} on $f$ is defined to be the subgroupoid of the homotopy fiber of:
	\[\Alg_{\bP^{\gr}_{[n+1,n]}}(\calC^{\cfil}) \rightarrow \left(\CAlg(\calC^{\cfil})\right)^{\Delta^1}\]
	at $f$.
\end{defn}
We then have a polyvector characterization of this space:
\begin{prop}\label{prop:coispoly}
	Let $f: A \rightarrow B$ a morphism in $\CAlg(\calC^{\cfil})$. Then we have an equivalence:
	\[\Cois^{\eps}(f,n) \simeq \Map_{\Alg_{\calL_{\infty}}(\Mod_k^{\cfil,\gr})}(k\langle 2 \rangle[-1],\wPol^{\eps}(f,n)\otimes k(-1)[n+1]).\]
	Moreover, it is compatible with the equivalence of \cref{prop:bivec}. The maps:
	\[\begin{tikzcd}
		&\Cois^{\eps}(f,n)\ar[rd]\ar[ld]&\\
		\Pois^{\eps}(A,n) && \Pois^{\eps}(B,n-1)
	\end{tikzcd}\]
	obtained from the morphisms of operads $\bP^{\gr}_{n+1} \rightarrow \bP^{\gr}_{[n+1,n]}$ and $\bP^{\gr}_n \rightarrow \bP^{\gr}_{[n+1,n]}$ are identified to the maps obtained by applying $\Map_{\Alg_{\calL_{\infty}}(\Mod_k^{\cfil,\gr})}(k\langle 2 \rangle[-1],\vari)$ to:
	\[\begin{tikzcd}
		&\wPol^{\eps}(f,n)^{\geq 2}[n+1]\ar[rd]\ar[ld]&\\
		\wPol^{\eps}(A,n)^{\geq 2}[n+1] && \wPol^{\eps}(B,n-1)^{\geq 2}[n].
	\end{tikzcd}\]
\end{prop}
\begin{proof}
	Thanks to all the setup of deformation theory of relative Swiss-Cheese algebras in the complete filtered context, the proof is an analog of \cite[4.16.]{MSI}.
\end{proof}
If $f: A \rightarrow B$ has a $k(-1)[n]$-twisted coisotropic structure and cotangent complexes of $A$ and $B$ are perfect, according to \cref{prop:coispoly} and the comparison between strict and weak relative polyvectors, we have natural maps:
\[\begin{tikzcd}
	\cot_A\otimes_A B(-1)[n] \ar[d] \ar[r] & \cot_B(-1)[n] \ar[r]\ar[d] & \cot_{B/A}(-1)[n]\ar[d]\\
	\tan_A\otimes_A B \ar[r] & \tan_{B/A}[1] \ar[r] & \tan_B[1].
\end{tikzcd}\]
We say that $f$ is \defterm{non-degenerate}, when $B$ has a perfect cotangent complex and $f$ has a perfect relative cotangent complex, if the underlying Poisson structure on $A$ is non-degenerate so that the first vertical map is an equivalence, and one of the two other vertical maps is an equivalence. We denote by $\Cois^{\eps,nd}(f,n)$ the maximal subgroupoid spanned by non-degenerate morphisms. The same comment on checking non-degeneracy on the associated graded that was mentioned for Poissons structures applies as well here.
\begin{thm}\label{thm:coisndlag}
	We have an equivalence:
	\[\Cois^{\eps,nd}(f,n)\simeq \Lag(f,k(-1)[n]).\]
\end{thm}
\begin{proof}
	The proof is an analog of \cite{MSII} or \cite{coisndpri}. We mimic it in \cref{sec:proofnd}.
\end{proof}

\section{Proof of the main theorem}\label{sec:thrdstep}
\par We provide the proof of \cref{thm:main}. In order to fully prove it, one has to establish a result relating shifted Poisson algebras in $\Mod_k$ to shifted Poisson algebras of weight $-1$ in $\Mod_k$, that is the point of \cref{ssec:pois2epspois}. We also need to prove that Poisson structures and Lagrangian thickenings satisfy Zariski descent, this will be done in \cref{ssec:desc}. We finally prove the theorem in \cref{ssec:proof}.
\subsection{From Poisson algebras to complete filtered Poisson algebras}\label{ssec:pois2epspois}
\par The operads $\bP_n$ and $\bP^{\gr}_n$ are related in the following way: the inclusion $(\vari)(0): \calC \rightarrow \calC^{\cfil}$ has a left adjoint $|\vari|^l$. With this functor, we get that $\left|\bP^{\gr}_n\right|^l\simeq \bP_n$. For filtration degree reasons, it is possible to see any $n$-shifted Poisson structure on $A \in \CAlg(\calC)$ as a $k(-1)[n]$-twisted Poisson structure on $A(0)$. We write this statement as follows:
\begin{prop}\label{prop:zeroPois1Pois}
	There is a fully faithful embedding
	\[\Alg_{\bP_n}(\calC) \hookrightarrow \Alg_{\bP^{\gr}_n}(\calC^{\cfil})\]
	with essential image spanned by those algebras with underlying object being constant in degree $0$. It is built in such a way that the following diagram:
	\[\begin{tikzcd}
		\Alg_{\bP_n}(\calC) \ar[r,hook]\ar[d] & \Alg_{\bP^{\gr}_n}(\calC^{\cfil}) \ar[d]\\
		\CAlg(\calC) \ar[r,"(\vari)(0)"] & \CAlg(\calC^{\cfil})
	\end{tikzcd}\]
	commutes.
\end{prop}
\begin{proof}
	Because of the maps $k(i) \rightarrow k(0)$ for $i\leq 0$, we naturally have a map of symmetric sequences:
	\[\Lie\{k(1)\} \rightarrow \Lie.\]
	If $[\vari,\vari]$ denotes the $\Lie$ bracket, $\mathrm{jac}$ the Jacobi relation, and $\mathrm{Free}$ the free operad, it is explicitly obtained from the commutative diagram:
	\[\begin{tikzcd}
		\mathrm{Free}(k(-2).\mathrm{jac})\ar[d]\ar[r] & \mathrm{Free}(k.\mathrm{jac})\ar[d]\\	
		\mathrm{Free}(k(-1).[\vari,\vari]) \ar[r] & \mathrm{Free}(k.[\vari,\vari])
	\end{tikzcd}\]
	by taking the operad quotients. The left-hand side recovers $\Lie\{k(1)\}$, and the right-hand side recovers $\Lie$. We get, in a similar way, a map of operads
	\[\bP^{\gr}_n \rightarrow \bP_n\]
	in $\calC^{\cfil}$. Thus we have a functor:
	\[\Alg_{\bP_n}(\calC^{\cfil}) \rightarrow \Alg_{\bP^{\gr}_n}(\calC^{\cfil}).\]
	commuting with the forgetful functors. We claim that it induces an equivalence
	\[\Alg_{\bP_n}(\calC^{\cfil})\times_{\calC^{\cfil}}\calC \simeq \Alg_{\bP^{\gr}_n}(\calC^{\cfil})\times_{\calC^{\cfil}}\calC.\]
	Indeed, the fiber of the right-hand side over $X \in \calC$ is given by
	\begin{align*}
		\Map_{\Op_{\calC^{\cfil}}}(\bP^{\gr}_n,\End_{X(0)})&\simeq\Map_{\Op_{\calC^{\cfil}}}(\bP^{\gr}_n,\End_X(0))\\
								   &\simeq\Map_{\Op_{\calC}}(|\bP^{\gr}_n|^l,\End_X)\\
								   &\simeq\Map_{\Op_{\calC}}(\bP_n,\End_X)
	\end{align*}
	which is the fiber over $X$ of the left-hand side. The first equivalence comes from the computation
	\begin{align*}
		F^i\Hom(X(0),Y(0)) & \simeq \Hom(X(0),Y(0)\otimes \mathbf{1}_{\calC}(-i))\\
				   &\simeq \Hom(X(0),Y(-i))\\
				   &\simeq \begin{cases} \Hom(X,Y) & \text{ if } i\leq 0,\\ 0 &\text{otherwise.}\end{cases}\\
				   &\simeq F^i(\Hom(X,Y)(0))
	\end{align*}
	for any $X,Y\in \calC$. The second one comes from the fact that $|\vari|^l$ is symmetric monoidal on complete filtered objects fully in non-positive degrees (\cite[1.8.]{MSI}).
	We conclude because $\Alg_{\bP_n}(\calC^{\cfil})\times_{\calC^{\cfil}} \calC \simeq \Alg_{\bP_n}(\calC)$, and the way the functor is built is compatible with underlying commutative structures.
\end{proof}

\subsection{Descent for Poisson structures and Lagrangian thickenings}\label{ssec:desc}
\par We establish two descent results needed to prove the main theorem. Assume $X$ is a derived scheme with a perfect cotangent complex, and consider $U_\bullet \rightarrow X$ a Zariski hypercover by affine derived scheme. The algebra of polyvectors $\wPol(X,n)$ refers to the one defined in \cite[3.1.]{CPTVV}. As a consequence of \cite[2.3.8.]{CPTVV} we know that its underlying graded algebra is given by
\[\Gamma(X,\Sym_X(\tan_X[-n](1))).\]
Because $\wPol(\vari,n)$ is functorial with respect to formally étale morphisms, we get a diagram $\wPol(U_\bullet,n)$. We can check that
\[\wPol(X,n) \simeq \lim \wPol(U_\bullet,n)\]
as graded $\Lie$-algebras. Indeed, as graded $\Lie$-algebras, we have a natural map
\[\wPol(X,n) \rightarrow \lim \wPol(U_\bullet,n)\]
to check it is an equivalence we can check it on each graded piece. But it is an equivalence on each piece since the tangent complex satisfies formally etale descent. Using \cite[3.1.2.]{CPTVV} and the fact that $\Map_{\Alg_{\Lie}(\calC^\gr)}(k\langle 2 \rangle[-1],\vari)$ preserves limits, we get:
\begin{prop}\label{prop:poisdesc}
	\[\Pois(X,n) \simeq \lim \Pois(U_\bullet,n).\]
\end{prop}
\par Now, assume that $X$ is locally of finite presentation and assume that $U_\bullet$ is composed only of affine derived schemes of finite presentation. Thanks to \cite[5.22.]{deffol}, we know that thickenings satisfy Zariski descent:
\[\Thick(X) \simeq \lim \Thick(U_\bullet).\]
To establish the result for Lagrangian structures, one has to prove that it is true for isotropic structures before, but this is a consequence of the fact that the Rham complex satisfies smooth descent:
\[\DR(X) \simeq \lim \DR(U_\bullet).\]
A proof of this result can be found, for instance, in the proof of \cite[1.14.]{PTVV} or \cite[2.6.]{sympgrpd}.
\begin{prop}\label{prop:lagthickdesc}
	\[\LagThick(X,n)\simeq\lim \LagThick(U_\bullet,n).\]
\end{prop}
\begin{proof}
	Let's unfold how the equivalence
	\[\Thick(X) \simeq \lim \Thick(U_\bullet)\]
	is constructed. The map $\Thick(X) \rightarrow \lim \Thick(U_\bullet)$ sends a thickening $X \rightarrow Y$ to the family $(Y^\wedge_U)$. Because $Y^\wedge_U \rightarrow Y$ is formally etale, it induces a map
	\[\LagThick(X,n) \rightarrow \lim\LagThick(U_\bullet,n).\]
	The converse map $\lim \Thick(U_\bullet) \rightarrow \Thick(X)$ sends a family $(Y_U)$ of thickenings to
	\[((\colim Y_U)^{\mathrm{Zar}})^{\conv}\]
	where $(\vari)^{\mathrm{Zar}}$ is the Zariski sheafification and $(\vari)^{\conv}$ is the convergent completion. By the universal property of the colimit, a compatible family of commutative squares
	\[\begin{tikzcd}
		U \ar[d] \ar[r] & Y_U \ar[d]\\
		* \ar[r,"0"] & \calA^{2,\cl}(\vari,n)
	\end{tikzcd}\]
	is the same datum as a commutative square
	\[\begin{tikzcd}
		X\simeq\colim U \ar[d] \ar[r] & \colim Y_U \ar[d]\\
		* \ar[r,"0"] & \calA^{2,\cl}(\vari,n).
	\end{tikzcd}\]
	It is then the same datum as commutative square
	\[\begin{tikzcd}
		X \ar[d] \ar[r] & (\colim Y_U)^{\mathrm{Zar}} \ar[d]\\
		* \ar[r,"0"] & \calA^{2,\cl}(\vari,n).
	\end{tikzcd}\]
	because $\calA^{2,\cl}(\vari,n)$ satisfies Zariski descent, and it is the same datum as
	\[\begin{tikzcd}
		X \ar[d] \ar[r] & ((\colim Y_U)^{\mathrm{Zar}})^{\conv} \ar[d]\\
		* \ar[r,"0"] & \calA^{2,\cl}(\vari,n).
	\end{tikzcd}\]
	essentially because $\DR(\vari)$ is convergent (see \cref{lm:DRconv}). This means that $\Isot\Thick(X,n)\simeq\lim \Isot\Thick(U_\bullet,n)$, and that
	\[\LagThick(X,n) \simeq \lim \LagThick(U_\bullet,n)\]
	since non-degeneracy of forms can be checked Zariski-locally.
\end{proof}

\subsection{Proof of the theorem}\label{ssec:proof}
\par We fix $A$ a commutative $k$-algebra with a perfect cotangent complex, we also fix $\tilde A$ to be a cofibrant strict model of $A$ (automatically bifibrant because every connective commutative differential graded $k$-algebra is fibrant).

\begin{lm}\label{lm:centernondegenerate}
	Assume $A(0)$ has a $\bP^{\gr}_{n+1}$-algebra structure. Then the natural $(n+1)$-shifted coisotropic morphism of weight $-1$ $\Z(\tilde A(0))\rightarrow \tilde A(0)$ is non-degenerate.
\end{lm}
\begin{proof}
	We have that, $\Z(\tilde A(0)) \simeq \rZ(\tilde A(0))$ is given by $\Pol^{\eps,\int}(A(0),n)^{\pi}$ where $\pi$ is the Maurer--Cartan element giving the Poisson structure of $A(0)$. The $\bP_{n+2}$-structure on $\rZ(\tilde A(0))$ is given by the Schouten--Nijenhuis bracket, but because non-degeneracy can be checked on the underlying graded object, and that the underlying object is given by usual graded internal polyvector, the center is non-degenerate. The coisotropic structure is also non-degenerate for the same reason.
\end{proof}

\begin{remark}
	The $(n+1)$-shifted coisotropic morphism $\rZ(A) \rightarrow A$ is in general not non-degenerate. In fact, even if $A$ has a perfect cotangent complex, $\rZ(A)$ does not need to have a perfect cotangent complex. For instance, if $A = k[x]$ with the trivial Poisson structure, $\rZ(A)$ is the formal power series ring $k[x][\![\partial_x]\!]$, with $\partial_x$ in degree $n+1$. When $n=-1$ it is a formal ring, and the cotagent complex is not perfect. Hence we have to take into account the formal structure of $\rZ(A)$, which is actually given by the complete filtered structure of $\rZ(A(0))$, where $A(0)$ is seen as a $\bP_{n+1}$-algebra of weight $-1$. It will actually give a derived foliation. One can then check that $|\rZ(A(0))| \simeq Z(A)$.
\end{remark}

\begin{thm}\label{thm:main}
	If $X$ is a locally of finite presentation derived scheme, there is an equivalence:
	\[\Pois(X,n) \simeq \LagThick(\Spec(X),n+1).\]
\end{thm}
\begin{proof}
	First, we prove it in the affine case for $A$ of finite presentation. We proceed in several steps. Let
	\[\Cois^{\eps}(A(0),n)\]
	be the space of maps $B \rightarrow A(0)$ of complete filtered commutative $k$-algebras together with an $n$-shifted coisotropic structure of weight $-1$. We let $\Cois^{\eps,nd}(A(0),n)$, $\Cois^{\eps,q.f.}(A(0),n)$ be the full subgroupoids spanned, respectively, by maps with non-degenerate structures and by maps $B\rightarrow A(0)$ such that $(\gr B)^0 \rightarrow A$ is an equivalence, and $\Sym_{(\gr B)^0)}((\gr B)^1(1))\rightarrow \gr B$ is an equivalence. The symbol $q.f.$ stands for ``quasi-free". Let
	\[\Lag^{\eps}(A(0),n)\]
	be the space of maps $B \rightarrow A(0)$ together with an $n$-shifted Lagrangian structure of weight $-1$. We set $\Lag^{\eps,q.f.}(A(0),n)$ to be the full subgroupoid defined in a similar way as $\Cois^{\eps,q.f.}(A(0),n)$. Now, by \cref{prop:zeroPois1Pois}, there is a commutative square:
	\[\begin{tikzcd}
		\Alg_{\bP_n}(\Mod_k) \ar[r,hook]\ar[d] & \Alg_{\bP^{\gr}_n}(\Mod_k^{\cfil}) \ar[d]\\
		\CAlg(\Mod_k) \ar[r,"(\vari)(0)"] & \CAlg(\Mod_k^{\cfil}).
	\end{tikzcd}\]
	By taking the fiber at $A$, we get the equivalence:
	\[\Pois(A,n) \simeq \Pois^{\eps}(A(0),n).\]
	Using the same technique as \cite[1.4.15.]{CPTVV}, we build a functor:
	\[\rZ: \Alg_{\bP^{\gr}_{n+1}}(\Mod_k^{\cfil})^{\mathrm{fet}} \rightarrow \Alg_{\bP^{\gr}_{[n+2,n+1]}}(\Mod_k^{\cfil})\]
	that is given by $B \mapsto (\rZ(B) \rightarrow B)$, when $B$ is cofibrant complete filtered $\bP^\gr_{n+1}$-algebra and where $\rZ(B)$ is the strict center. The superscript $\mathrm{fet}$ stands for ``formally etale" and means that we restrict to the category with only formally etale morphisms. It is in fact constructed as follows:
	\begin{itemize}
		\item Set $I$ to be the diagram of all complete filtered $\bP^\gr_{n+1}$-algebras with formally etale morphisms;
		\item The forgetful functor $I \rightarrow \strMod_k^{\cfil}$ make this diagram into a $\bP^\gr_{n+1}$-algebra in $\mathrm{Fun}(I,\strMod_k^{\cfil})$ which is a model category coming from the injective model structure (there are some size issues here but, for our purpose, we can fix the algebras to be of cardinality at most the cardinality of $A$);
		\item We pick a fibrant-cofibrant resolution of this algebra that we will denote by $\underline{B}$, and take the (internal) strict center $\rZ(\underline{B})$;
		\item This models our functor $\rZ$ and is the correct one thanks to \cite[1.4.13.]{CPTVV}.
	\end{itemize}
	It fits into a commutative square:
	\[\begin{tikzcd}
		\Alg_{\bP^{\gr}_{n+1}}(\Mod_k^{\cfil})^{\mathrm{fet}} \ar[rr] \ar[rd] && \Alg_{\bP^{\gr}_{[n+2,n+1]}}(\Mod_k^{\cfil}) \ar[ld]\\
		&\Alg_{\Com}(\Mod_k^{\cfil})^{fet}
	\end{tikzcd}\]

	By taking the fiber at $A(0)$, we get a map
	\[\rZ : \Pois^{\eps}(A(0),n) \rightarrow \Cois^{\eps}(A(0),n+1)\]
	However, because the center is quasi-free and the coisotropic structure is non-degenerate according to \cref{lm:centernondegenerate}, it restricts to a map:
	\[\Pois^{\eps}(A(0),n) \rightarrow \Cois^{\eps,nd,q.f.}(A(0),n+1).\]
	Now, there is a converse map:
	\[P: \Cois^{\eps,nd,q.f.}(A(0),n+1) \rightarrow \Pois^{\eps}(A(0),n)\]
	which is the forgetful functor to Poisson structures on the target given in \cref{prop:coispoly}. It is clear that $P \circ \rZ \simeq \id$. For the converse, given a coisotropic morphism $B \rightarrow A(0)$, one has naturally a map $z: B \rightarrow \rZ(P(A(0)))$ by definition of coisotropic structures. Now, because the morphism is non-degenerate, we have that $\gr(B)^1\simeq \tan_A[-n+1]$. So $z$ is an equivalence on the associated graded and then an equivalence. We now have
	\[\Pois^{\eps}(A(0),n)\simeq \Cois^{\eps,nd,q.f.}(A(0),n+1).\]
	By \cref{thm:coisndlag}, there is an equivalence
	\[\Lag^{\eps}(A(0),n+1)\simeq\Cois^{\eps,nd}(A(0),n+1)\]
	which restricts into an equivalence
	\[\Lag^{\eps,q.f.}(A(0),n+1)\simeq\Cois^{\eps,nd,q.f.}(A(0),n)\]
	There is also an equivalence 
	\[\Lag^{\eps,q.f.}(A(0),n+1) \simeq \LagThick(\Spec(A),n+1)\]
	by \cref{thm:lagalglagthick}. To conclude, we apply \cref{prop:poisdesc} and \cref{prop:lagthickdesc}: if $U_\bullet \rightarrow X$ is a Zariski hypercover of $X$ composed of affine derived schemes of finite presentation, then
	\[\Pois(X,n) \simeq \lim \Pois(U_\bullet,n) \simeq \lim \LagThick(U_\bullet,n+1)\simeq \LagThick(X,n+1).\]
\end{proof}
\par As a corollary, we can prove the following:
\begin{cor}\label{cor:aksz}
	Let $Y$ be a $\calO$-compact with $\calO$-orientation of degree $d$ derived stack, in the sense of \cite[2.1.]{PTVV}. Let $X$ be a locally of finite presentation derived scheme with an $n$-shifted Poisson structure. Then there exists a map $\Map(Y,X) \rightarrow Z$ with an $(n-d+1)$-shifted Lagrangian structure.	If $Y = M_{B}$ is the Betti stack of a compact oriented $d$-dimensional manifold, then there exists an $(n-d)$-shifted Poisson structure on $\Map(M_B,X)$.
\end{cor}
\begin{proof}
	By \cref{thm:main}, we can consider a Lagrangian thickening $X \rightarrow X^\mathrm{symp}$. It has an $(n+1)$-shifted Lagrangian structure and then
	\[f : \Map(Y,X) \rightarrow \Map(Y,X^\mathrm{symp})\]
	has a natural $(n+1-d)$-shifted Lagrangian structure by \cite[2.35. (2)]{sympgrpd}. Now assume $Y = M_B$ is a Betti stack of a compact oriented manifold of dimension $d$. Since $M$ is finite, $\Map(M_B,X)$ is a finite limit of a constant diagram with value $X$. Then it is a derived scheme locally of finite presentation and $f$ is locally almost of finite presentation. It is also a formal thickening:
	\begin{align*}
		\Map(M_B,X)_\dR &\simeq (\lim_M X)_\dR\\
				&\simeq \lim_M (X_\dR)\\
				&\simeq \lim_M (X^\mathrm{symp}_\dR)\\
				&\simeq \Map(M_B,X^\mathrm{symp})_\dR
	\end{align*}
	and $\Map(M_B,X^\mathrm{symp})$ has a deformation theory for the same reason. The map $f$ is a Lagrangian thickening; we then apply \cref{thm:main} to extract a Poisson structure on $\Map(M_B,X)$.
\end{proof}
\begin{cor}
	Let $L \rightarrow X$ be a $n$-shifted Lagrangian morphism, then the formal neighborhood of $L$ is equivalent to the zero section of
	\[L \rightarrow \widehat{T^*[n](L/X)}\]
	if and only if the underlying $(n-1)$-shifted Poisson structure on $L$ is trivial.
\end{cor}
\begin{example}
	Fix $t \in \C$. There is a $1$-parameter family of actions of $\Z^2$ on $\C\times\widehat\C$ given by
	\[(m,n).(z,w) \defeq (z + m + (i+t w)n,w).\]
	Take the torus $T \defeq \C/\Z^2$ with the action above at $t = 0$. We also set $Y_t \defeq (\C\times\widehat\C)/\Z^2$ with the symplectic form $\ddr z \wedge \ddr w$. We then get a family of Lagrangian submanifolds
	\[T \hookrightarrow Y_t.\]
	The Lagrangian structure gives a $(-1)$-shifted Poisson structure on $T$. We can check for degree reasons that such structures are actually parametrized by $H^1(T,\Sym^2(\tan_T)) \simeq H^1(T,\calO_T)$ because $T$ is an algebraic group of dimension $1$, which is $1$-dimensional by Riemann--Roch. When $t = 0$, $Y_0$ is actually the (formal) cotangent bundle of $T$ and $T \rightarrow Y_0$ is the zero section, the associated $(-1)$-shifted Poisson structure on $T$ is the trivial one. When $t \neq 0$, the $(-1)$-shifted Poisson structure on $T$ is not trivial and $Y_t$ can't be equivalent to the cotangent bundle of $T$.
\end{example}
\section{Future developments}
\par The technique of formal localization is well-behaved for Artin derived stacks. Using the tools developed by Brantner, Magidson, and Nuiten (\cite{deffol}), we claim that it is possible to generalize theorems of formal localization to locally of finite presentation prestacks with deformation theories and perfect cotangent complexes. We plan to investigate this in the future and generalize to any well-behaved prestacks the main theorem of this article. This will allow us to show the existence of new Poisson structures, specifically from AKSZ-like constructions coming from non-necessarily compact stacks. For instance, Pantev and To\"en (\cite{akszncompdr, akszncompbet}) build a shifted Lagrangian morphism between stacks that are not Artin but still have a well-behaved deformation theories.
\par In the relative setting, one may also look for a similar theorem for shifted coisotropic structures. This question is related to shifted coisotropic reductions, generalizing shifted Hamiltonian reductions built by Anel and Calaque (\cite{hamred}). We can conjecture the following:
\begin{conjecture}
	Let $f: X \rightarrow Y$ be a finitely presented morphism of affine derived $k$-schemes, with $Y$ of finite presentation. Then the space of $n$-shifted coisotropic structures on $f$ is equivalent to the space of diagrams:
	\[\begin{tikzcd}
	X\\
	& X^{\mathrm{symp}} & Y \\
	& X//f & Y^{\mathrm{symp}}
	\arrow[from=1-1, to=2-2, "\pi_X"]
	\arrow[bend left, from=1-1, to=2-3, "f"]
	\arrow[bend right, from=1-1, to=3-2, "\red"]
	\arrow[from=2-2, to=2-3]
	\arrow[from=2-2, to=3-2]
	\arrow["\lrcorner"{anchor=center, pos=0.125}, draw=none, from=2-2, to=3-3]
	\arrow[from=2-3, to=3-3, "\pi_Y"]
	\arrow[from=3-2, to=3-3, "\pi_f"]
	\end{tikzcd}\]
	where:
	\begin{itemize}
		\item $\pi_Y$ is an $(n+1)$-shifted Lagrangian thickening,
		\item $\pi_f$ is an $(n+1)$-shifted Lagrangian morphism,
		\item $\pi_X$ is an $n$-shifted Lagrangian thickening,
		\item $\red$ is a thickening, and
		\item we have an identification between the two $n$-shifted symplectic structures on $X^{\mathrm{symp}}$. One is induced from the Lagrangian structure on $\pi_X$ and the other one is induced from the Lagrangian intersection.
	\end{itemize}
	We call the prestack $X//f$ the \defterm{coisotropic reduction} of $f$.
\end{conjecture}
\begin{proof}{(Sketch)}
	We sketch here a strategy to construct such a diagram. First, set $X=\Spec(A)$ and $Y=\Spec(B)$. Then $B\rightarrow A$ has a $\bP_{[n+1,n]}$-algebra structure. Naturally, we then have a map $B \rightarrow \Z(A)$ and we know that we have a weak mixed structure on $\wPol(A/B,n-2)$. We denote by $\Z_B(A)$ this weak graded algebra that may be identified with a complete filtered object. We should have the following diagram:
	\[\begin{tikzcd}
		A \\
		& \Z(A) & B \\
		& \Z_B(A) & \Z(B)
		\arrow[from=2-2, to=1-1, "p_A"']
		\arrow["\lrcorner"{anchor=center, pos=0.125}, draw=none, from=2-2, to=3-3]
		\arrow[from=2-3, to=1-1, bend right, "f^\sharp"']
		\arrow[from=2-3, to=2-2]
		\arrow[from=3-2, to=1-1, bend left]
		\arrow[from=3-2, to=2-2	]
		\arrow[from=3-3, to=2-3, "p_f"']
		\arrow[from=3-3, to=3-2, "p_B"']
	\end{tikzcd}\]
	in complete filtered objects. We expect that $p_f$ has a non-degenerate $(n+1)$-shifted coisotropic structure of weight $-1$ and that the induced coisotropic structure on $\Z(A)$ coming from coisotropic intersection is the same as the one coming from the Schouten--Nijenhuis bracket. Once this is settled, one has to check that the non-degenerate coisotropic intersection is the same thing as the Lagrangian intersection, which is a conjecture stated in \cite[3.4.4.]{coiscorr}. Because $\Z_B(A)$ is a derived foliation over $A$, we can then apply a theorem in the same fashion as \cref{thm:lagalglagthick} to induce the diagram of the conjecture.
\end{proof}
In the case $Y$ is shifted symplectic, this conjecture was originally announced by Costello--Rozenblyum\footnote[1]{\textit{BV formalism and derived symplectic geometry}, \url{https://doi.org/10.48660/17010004}.}.
The author also believes that for critical loci, it is possible to relate some of the construction involved here to the Classical Master Equation of BV formalism.
\appendix
\section{Non-degenerate Poisson and coisotropic structures are symplectic and Lagrangian}\label{sec:proofnd}
\par This section ought to sketch in the complete filtered context the proof of \cref{thm:poisndsymp} and \cref{thm:coisndlag}. Let's give a general lemma before:
\begin{lm}\label{lm:mcepstw}
	If $\fg$ is a $\Lie\{k\langle 1 \rangle(1)\}$-algebra in $\Mod_k^{\gr,\cfil}$, then $\fg\otimes k[\eta]$, where $\eta$ is of degree $0$ and satisfies $\eta^2 = 0$, is also a $Lie\{k\langle 1 \rangle(1)\}$-algebra. For convenience, we denote by $\Lie$ the category $\Alg_{\Lie\{k\langle 1 \rangle(1)\}}(\Mod_k^{\gr,\cfil})$. There is a pullback diagram:
	\[\begin{tikzcd}
		\Map_{\Mod_k^{\eps,\gr,\cfil}}(k\langle 2 \rangle [-1],\fg_x) \ar[r] \ar[d] & \Map_\Lie(k\langle 2 \rangle[-1],\fg\otimes k[\eta])\ar[d]\\
		* \ar[r,"x"] & \Map_\Lie(k\langle 2 \rangle[-1],\fg)
	\end{tikzcd}\]
	where $\fg_x$ is the $\Lie$-algebra $\fg$ but with an extra mixed structure given by $[x,-]$ (The upper-left corner is computed by hand by choosing the projective model structure to describe $\Mod_k^{\eps,\gr,\cfil}$ and consider a cofibrant replacement of $k\langle 2 \rangle[-1]$, instead of taking the injective model structure).
\end{lm}
The proof is exactly the same as in the unfiltered case; we can compute the Maurer--Cartan element defined by $x$ that is in filtration degree $0$ and twist as usual without worrying of the extra filtration. Because of the following proposition, it will be enough to prove \cref{thm:coisndlag} to prove \cref{thm:poisndsymp}. Though it is also possible to prove it in a more direct way by mimicking the proof of \cite[3.2.5.]{CPTVV} or \cite[1.38.]{poisndpri}. For convenience, we set $\Lag^{\eps}(f,n) \defeq \Lag(f,k(-1)[n])$.
\begin{prop}\label{prop:cois2pois}
	If for any $f: A\rightarrow B$, $\Cois^{\eps,nd}(f,n) \simeq \Lag^{\eps}(f,n)$ then for any $A$ $\Pois^{\eps,nd} \simeq \Symp(A,n)$.
\end{prop}
\begin{proof}
	It's because we can apply the first equivalence to $f: k \rightarrow A$. It is straightforward to see that $\Lag^{\eps}(f,n-1) \simeq \Symp(A,n)$. There is also an equivalence $\Cois^{\eps}(f,n-1) \simeq \Pois^{\eps}(A,n)$. Indeed, a coisotropic structure on $f$ is the data of a $\bP^{\gr}_{[n+1,n]}$-algebra structure on $f$, that is a $\bP^{\gr}_n$-algebra structure on $k$, a $\bP^{\gr}_{n+1}$-algebra on $A$ structure and a map $\Z(k) \rightarrow A$ extending $f$. But $\Z(k) \simeq k$ so this gives no extra data. We also see directly that by definition of non-degeneracy, $\Cois^{\eps,nd}(f,n-1) \simeq \Pois^{\eps,nd}(A,n)$.
\end{proof}
We then recall the main ingredients of \cite[4.22.]{MSII} and \cite[3.17.]{coisndpri}.

\par Assume $f: A\rightarrow B$ is a morphism between cofibrant objects in the model category $\strCAlg(\rmM^{\cfil})$. We know that there is a homotopy graded $\bP^{\gr}_{[n+2,n+1]}$-algebra consisting of weak polyvectors
\[(\wPol^{\eps}(A,n),\wPol^{\eps}(B/A,n-1)).\]
We then do as in \cite{MSII} and pick the strict model of strict polyvectors, which is a strict graded $\bP^{\gr}_{[n+2,n+1]}$-algebra:
\[(\Pol^{\eps}(A,n),\Pol^{\eps}(B/A,n-1)).\]
We let $\Pol^{\eps}(f,n)$ be, as a graded commutative algebra, the homotopy fiber of $\Pol^{\eps}(A,n) \rightarrow \Pol^{\eps}(B/A,n)$. Given a shifted coisotropic structure of weight $-1$ on $f$, that is a Maurer--Cartan element $(\gamma_A,\gamma_B)$ of $\Pol^{\eps}(f,n)$ (if we consider moreover its $\calL_{\infty}$-structure), we get a diagram:
\[\begin{tikzcd}
	{|\DR(A)|^{\cfil}} \ar[r] \ar[d,"\mu(\vari{,}\gamma)_A"] & {|\DR(B)|^{\cfil}} \ar[d,"\mu(\vari{,}\gamma)_B"]\\
	\Pol^{\eps}(A,n)^{\gamma_A} \ar[r] &\Pol^{\eps}(B/A,n-1)^{\gamma_B}.
\end{tikzcd}\]
Both polyvectors have weak mixed structure given by $[\gamma_A,\vari]$ for the first one and $[\gamma_B,\vari] + \sum_{k\geq 1}\frac1{(k-1)!}f_k(\gamma_A;\gamma_B,\dots,\gamma_B,\vari)$, where the $f_k$ comes from the map to the center. This has been already described in \cite[3.6.]{MSI}, but we only have an extra filtration, which is ok since everything is compatible with the filtration. It is then easy to check that $f$ is non-degenerate if and only if the vertical arrows are equivalences. If we consider the fibers we have a map $\mu(\vari,\gamma): \DR(f) \rightarrow \Pol(f,n)^{\gamma}$ where $\DR(f) = \fib(|\DR(A)|^{\cfil} \rightarrow |\DR(B)|^{\cfil})$ as a graded mixed object. Now, we define the space $\TCois^{\eps}(A,n)$:
\[\TCois^{\eps}(A,n) \defeq \Map_{\Alg_{\Lie\{k\langle 1\rangle(1)\}}(\Mod_k^{\gr,\cfil})}(k\langle 2 \rangle[-1],\Pol^{\eps}(f,n)\otimes k(-1)[n+1]\otimes k[\eta]).\]
\par For simplicity, we shall omit the parameters $(A,n)$ in each notation and leave it implicit. The projection:
\[p: \TCois^{\eps} \times \Lag^{\eps} \rightarrow \Cois^{\eps} \times \Lag^{\eps}\]
has a section given by
\[\Phi: (\gamma,\lambda) \mapsto (\gamma + \eta(\sigma \gamma - \mu(\lambda,\gamma)),\lambda)\]
where $\sigma$ sends an element of degree $p$ to $(p-1)$ times this element. For the same reason as in \cite[4.17.]{MSII}, it is well defined, and we can define the space of compatible pairs $\Comp^{\eps}(f,n)$ to be the relative zero locus of this section, that is, the homotopy limit of:
\[\Cois^{\eps} \times \Lag^{\eps} \overset{\Phi}{\underset{0}\rightrightarrows} \TCois^{\eps} \times \Lag^{\eps} \overset{p}\rightarrow \Cois^{\eps}\times\Lag^{\eps}.\]
We set $\Comp^{\eps,nd}$ to be the space of compatible pairs with underlying coisotropic structure non-degenerate. We have the following statement:
\begin{prop}\label{prop:comp2cois}
	The natural map $\Comp^{\eps,nd} \rightarrow \Cois^{\eps,nd}$ is an equivalence.
\end{prop}
This is because, by definition of non-degenerate structures, if $\gamma$ is non-degenerate; then $\mu(\vari,\gamma): \DR(f) \rightarrow \Pol(f,n)^{\gamma}$ is an equivalence. The hard statement is the other one requiring an induction:
\begin{prop}\label{prop:comp2lag}
	The natural map $\Comp^{\eps,nd} \rightarrow \Lag^{\eps}$ is an equivalence.
\end{prop}
We proceed by induction. Indeed, all the constructions we have done come with extra filtrations coming from the mixed graded structures (and not the internal filtration). $(\Lag^{\eps,\leq p} \subset \Map(k(1)\langle 2 \rangle[-n-2],\DR^{\leq p}(f)))$ for $\Lag^{\eps}$, $(\Cois^{\eps,\leq p})$ for $\Cois^{\eps}$, and  $(\Comp^{\eps,\leq p})$ for $\Comp^{\eps}$. It suffices then to show that $\Comp^{\eps,nd,\leq p} \rightarrow \Lag^{\eps,\leq p}$ is an equivalence for any $p \geq 2$ to conclude.
\begin{prop}\label{prop:init}
	\[\Comp^{\eps,nd,\leq 2} \rightarrow \Lag^{\eps,\leq 2}\]
	is an equivalence.
\end{prop}
\begin{proof}
	If $(\gamma,\lambda) \in \Comp^{\eps}$, it is possible to describe the weight $2$ part, which is the data of maps
	\[\gamma_A^{\sharp}: \cot_A(-1)[n] \rightarrow \tan_A, \qquad\gamma_B^{\sharp}: \cot_B(-1)[n-1] \rightarrow \tan_{B/A},\]
	and
	\[\lambda_A^{\sharp}: \tan_A \rightarrow \cot_A(-1)[n],\qquad \lambda_B^{\sharp}: \tan_{B/A} \rightarrow \cot_B(-1)[n-1].\]
	with equivalences coming from the compatibilities:
	\[\lambda_A^{\sharp} \circ \gamma_A^{\sharp} \circ \lambda_A^{\sharp} \simeq \lambda_A^{\sharp},\]
	and
	\[\lambda_B^{\sharp} \circ \gamma_B^{\sharp} \circ \lambda_B^{\sharp} \simeq \lambda_B^{\sharp}.\]
	To see that, we do exactly the same computation as in \cite[4.25.]{MSII}. Because $\lambda_A^{\sharp}$ and $\lambda_B^{\sharp}$ are equivalences, $\gamma_A^{\sharp}$ and $\gamma_B^{\sharp}$ are completely determined by those and give the weight $2$ part of a non-degenerate shifted coisotropic structure of weight $-1$.
\end{proof}
\par We now introduce obstruction spaces to lift.
\[\Obs(p+1,\Lag^\eps) \defeq F^0\DR(f)^{p+1}[n+3](-1)\times\Lag^{\eps,\leq p},\]
\[\Obs(p+1,\Cois^\eps) \defeq F^0\Pol^{\eps}(f,n)^{p+1}[n+3](-1)\times\Cois^{\eps,\leq p}\text{, and}\]
\[\Obs(p+1,\TCois^\eps) \defeq F^0\Pol^{\eps}(f,n)^{p+1}[n+3](-1)\otimes k[\eta]\times\TCois^{\eps,\leq p}.\]
There is a section
\[s: \Lag^{\eps,\leq p} \rightarrow \Obs(p+1,\Lag^\eps)\]
sending $\lambda = \Sigma_{i\leq p} \lambda_i$ to $(\ddr\lambda_p,\lambda)$. From this perspective, it is possible to extend an element $\lambda$ of $\Lag^{\eps,\leq p}$ to an element of $\Lag^{\eps,\leq p+1}$ if and only if $\ddr \lambda_p \simeq 0$.
\begin{lm}
	$\Lag^{\eps,\leq p+1}$ is the relative zero locus of
	\[\Obs(p+1,\Lag^\eps) \overset{s}\leftrightarrows \Lag^{\eps,\leq p}.\]
\end{lm}
Similarly for $\Cois^{\eps,\leq p}$ and $\TCois^{\eps,\leq p}$, one defines the section $\gamma \mapsto ((\frac12 P_{p+1}\{\gamma,\gamma\},\gamma)$, where $P_{p+1}$ is the weight $p+1$ part. The relative zero locus of these sections gives, respectively, $\Cois^{\eps,\leq p+1}$ and $\TCois^{\eps,\leq p+1}$. If we combine everything, we get the diagram:
\[\begin{tikzcd}
	\Obs(p+1,\TCois^\eps\times\Lag^\eps)\ar[r]\ar[d] & \Obs(p+1,\Cois^\eps\times\Lag^\eps)\ar[d]\\
	\TCois^{\eps,\leq p}\times\Lag^{\eps,\leq p} \ar[r] \ar[u,bend left] & \Cois^{\eps,\leq p}\times\Lag^{\eps,\leq p} \ar[l,bend right, "\Phi_p"] \ar[u,bend right]\\
	\TCois^{\eps,\leq p+1} \times \Lag^{\eps,\leq p+1} \ar[r]\ar[u] & \Cois^{\eps,\leq p+1} \times \Lag^{\eps,\leq p+1}. \ar[l,bend right, "\Phi_{p+1}"]\ar[u]
\end{tikzcd}\]
We build a section of the top horizontal arrow:
\[(\delta_\gamma,\delta_\lambda,\gamma,\lambda)\mapsto(\delta_\gamma + \eta(\sigma\delta_\gamma-\mu(\delta_\lambda,\gamma)-\xi(\lambda,\gamma,\delta_\gamma)),\delta_\lambda,\gamma + \eta(\sigma\gamma-\mu(\lambda,\gamma)),\lambda)\]
where $\xi(\vari,\gamma,\delta_\gamma)$ is the walking derivation relative to $\mu(\vari,\gamma)$ and given in weight $1$ by $\xi(a\ddr x,\gamma,\delta_\gamma)\defeq a[\delta_\gamma,x]$. It is clear that it is a section. By \cite[4.30.]{MSII} it is compatible with other sections. We define $\Obs(p+1,\Comp^{\eps})$ to be the relative zero locus of this section. Because it is compatible with other sections and by commutation of the limits, we get a relative zero locus diagram
\[\Obs(p+1,\Comp^\eps) \leftrightarrows \Comp^{\eps,\leq p} \leftarrow \Comp^{\eps,\leq p+1}.\]
We set $\Obs(p+1,\Comp^{\eps,nd})$ be the fiber of $\Obs(p+1,\Comp^\eps)$ over $\Comp^{\eps,nd,\leq p}$.
\begin{prop}\label{prop:ind}
	The projection
	\[\Obs(p+1,\Comp^{\eps,nd})\rightarrow \Obs(p+1,\Lag^\eps)\]
	is an equivalence.
\end{prop}
\begin{proof}
	The fiber of $\Obs(p+1,\Comp^{\eps})$ over $(\gamma,\lambda)\in \Comp^{\eps,\leq p}$ is given by the fiber of
	\[\begin{array}{ccc}
		\Pol^{\eps}(f,n)^{p+1}\oplus\DR(f)^{p+1} & \rightarrow & \Pol^{\eps}(f,n)^{p+1}\\
		(\gamma_{p+1},\lambda_{p+1})&\mapsto& p\gamma_{p+1}-\xi(\lambda_2,\gamma_2,\gamma_{p+1}) - \mu(\lambda_{p+1},\gamma_2).
	\end{array}\]
	In order to show the equivalence, it is enough to show that the fiber over a non-degenerate pair is $\DR(f)^{p+1}$. If we show that
	\[\begin{array}{ccc}
		\Pol^\eps(f,n)^{p+1} & \rightarrow & \Pol^\eps(f,n)^{p+1}\\
		x & \mapsto & p x - \xi(\lambda_2,\gamma_2,x)
	\end{array}\]
	is an equivalence then we see that $\gamma_{p+1}$ is completely determined by $\lambda_{p+1}$ and $\gamma_2$ and then it concludes the proof. We manage to do it the same way as in \cite[4.31.]{MSII}.
\end{proof}
\begin{proof}[Proof of \cref{prop:comp2lag}]
\par We have done the induction steps. The initial case is done in \cref{prop:init} and the induction in \cref{prop:ind}. It concludes because
\[\Comp^{\eps,nd} \simeq \lim_p \Comp^{\eps,nd,\leq p},\]
and
\[\Lag^{\eps} \simeq \lim_p \Lag^{\eps,\leq p}.\]
\end{proof}
\par We combine \cref{prop:comp2cois} and \cref{prop:comp2lag} to set once and for all \cref{thm:coisndlag}:
\[\Cois^{\eps,nd} \simeq \Lag^{\eps}.\]

\begin{acknowledgements}
	This work is part of a broader project aiming to understand BV formalism in derived geometry. I firstly warmly thank my advisor, Damien Calaque, for listening, helping, and trusting me and my ideas from the beginning. The main strategy to prove the theorem was given to me by Damien Calaque and Pavel Safronov. I am really grateful to Pavel Safronov for accepting to pursue his idea and giving me good insights. I would like to thank Joost Nuiten for his explanations about formal integration of Lie algebroids, pro-coherent modules, and for solving an issue I had about the existence of cotangent complex of formal thickenings. I am also grateful to Victor Alfieri, Mathieu Anel, Carlo Buccisano, Jiaqi Fu, Valerio Melani, Ulysse Mounoud, Hugo Pourcelot, Bertrand To\"en, and Sofian Tur-Dorvault for useful conversations related to this paper, and Alyosha Latyntsev for some remarks on a preliminary version of this paper. Finally, I would like to thank the anonymous referee for useful comments on a first version of this paper.
\end{acknowledgements}

\bibliographystyle{amsalpha}
\bibliography{bibliography}

\end{document}